\documentclass[11pt]{article}
\usepackage{amsmath}
\usepackage{amssymb}
\usepackage{amscd}
\usepackage{amsthm}
\usepackage{graphics}

\def\endproof{\relax\ifmmode\expandafter\endproofmath\else
  \unskip\nobreak\hfil\penalty50\hskip.75em\hbox{}\nobreak\hfil\bull
  {\parfillskip=0pt \finalhyphendemerits=0 \bigbreak}\fi}
\def\endproofmath$${\eqno\bull$$\bigbreak}
\def\bull{\vbox{\hrule\hbox{\vrule\kern3pt\vbox{\kern6pt}\kern3pt\vrule}\hrule}}
\newtheorem{theorem}{Theorem}[subsection]
\newtheorem*{main}{Theorem}
\newtheorem{proposition}[theorem]{Proposition}
\newtheorem{lemma}[theorem]{Lemma}

\newtheorem{corollary}[theorem]{Corollary}

\newtheorem{D}[theorem]{Definition}
\newenvironment{defn}{\begin{D} \rm }{\end{D}}

\newtheorem{R}[theorem]{Remark}
\newenvironment{remark}{\begin{R}\rm }{\end{R}}

\def\Zee{\mathbb{Z}}
\def\Q{\mathbb{Q}}
\def\Ar{\mathbb{R}}
\def\Cee{\mathbb{C}}

\def\scrO{\mathcal{O}}
\def\ov{\overline}
\def\spcheck{^{\vee}}
\def\frak{\mathfrak}

\def\Ker{\operatorname{Ker}}

\def\ad{\operatorname{ad}}

\def\Lie{\operatorname{Lie}}

\title{On the  Converse to a Theorem of Atiyah and Bott}
\author{Robert Friedman\thanks{The first author was partially
    supported by NSF grant DMS-99-70437.}  
\  and John W. Morgan\thanks{The second author was partially supported
   by NSF grant DMS-97-04507.}}

\begin{document}

\maketitle

\section*{Introduction}

Throughout this paper, $C$ denotes a smooth projective curve of genus at least
one, $G$ denotes a reductive linear algebraic group over $\Cee$, and $\xi_0$
is a $C^\infty$ principal  $G$-bundle
over $C$. The space of
all $(0,1)$-connections on $\xi_0$  is  an affine
space $\mathcal{A}=\mathcal{A}(\xi_0)$ associated to the infinite dimensional
complex vector space
$H^{0,1}(C;\ad \xi_0)$.
Following Shatz \cite{Shatz} for the case $G=GL(n)$,
Atiyah and Bott \cite{AtBo} defined a natural stratification
of this space. If $\frak h$ is a Cartan subalgebra for $G$,
then the strata are indexed by the orbits under the Weyl group of a certain
discrete set of points in $\frak h_\Ar$, the split real form of $\frak h$,
which we call points of {\sl
  Atiyah-Bott
type\/} for $\xi_0$. Fix a set $\Delta$ of simple roots for $G$ with respect
to
$\frak h$. Since every Weyl orbit in
$\frak h_\Ar$ has a unique representative in the positive Weyl chamber
$\ov C_0$ associated to $\Delta$, it is natural to index the strata by points
$\mu$ of Atiyah-Bott type for $\xi_0$ which lie in $\ov C_0$. We denote by
$\mathcal{C}_\mu$ the stratum corresponding to $\mu$. A point $\mu$ of
Atiyah-Bott type for $\xi_0$
 determines a parabolic subgroup $P(\mu)$ of $G$, together with a
$C^\infty$ $P(\mu)$-bundle $\eta_0(\mu)$ over $P(\mu)$ such that
$\eta_0(\mu)\times _{P(\mu)}G$ is $C^\infty$ isomorphic to $\xi_0$. The
condition that $\mu\in \ov C_0$ is just the condition that $P(\mu)$ is a
standard parabolic subgroup. Recall that an unstable holomorphic bundle $\xi$
whose underlying $C^\infty$-bundle is $C^\infty$ isomorphic to $\xi_0$ has a
canonical reduction to a parabolic subgroup, called the {\sl
Harder-Narasimhan\/} reduction. A
$(0,1)$-connection lies in
$\mathcal{C}_\mu$ if and only if the parabolic subgroup of its
Harder-Narasimhan reduction is conjugate to
$P(\mu)$ in such a way that the corresponding holomorphic $P(\mu)$-bundle is
$C^\infty$ isomorphic to $\eta_0(\mu)$. The strata $\mathcal{C}_\mu$ are of
finite codimension in $\mathcal{A}$  and are invariant under the action
of the group of $C^\infty$-changes of gauge.  We define
$\mathcal{C}_{\mu'}\preceq _r \mathcal{C}_{\mu}$ if the closure of
$\mathcal{C}_{\mu'}$ meets
$\mathcal{C}_\mu$, and define the relation
$\preceq$ on the strata
$\mathcal{C}_\mu$
by taking the unique extension of $\preceq _r$ to a transitive relation. More
concretely, $\mathcal{C}_{\mu'}\preceq   \mathcal{C}_{\mu}$ if and only if
there is a sequence
$\mathcal{C}_{\mu'}=\mathcal{C}_{\mu_0},\mathcal{C}_{\mu_1},
\ldots,\mathcal{C}_{\mu_n}=\mathcal{C}_{\mu}$ such that for each
$i\,,0\le i\le n-1$ the closure of $\mathcal{C}_{\mu_i}$ meets
$\mathcal{C}_{\mu_{i+1}}$.

For Atiyah and Bott, the stratification
arises as follows. Let $K$ be a compact Lie group whose complexification is
$G$, and let $\xi_K$ be a $C^\infty$ principal $K$-bundle such that
$\xi_K\times_KG$ is $C^\infty$ isomorphic to $\xi_0$. One can then identify
$K$-connections on $\xi_K$ with $(0,1)$-connections on $\xi_0$. Under this
identification, the strata $\mathcal{C}_\mu$  are the stable sets under
the negative gradient flow for the Yang-Mills functional $\|F_A\|^2$ on the
space of
$K$-connections on $\xi_K$ for the various connected components of the
critical set of the Yang-Mills functional.  From this point of view,
$\mathcal{C}_{\mu'}\preceq _r
\mathcal{C}_{\mu}$ if there is a downward gradient flow line in the stratum
$\mathcal{C}_{\mu'}$ converging as  $t \to -\infty$ to a point in
$\mathcal{C}_\mu$.

Atiyah and Bott introduce a partial ordering, the {\sl Atiyah-Bott
ordering\/}, on
  points $\mu\in \frak h$, as follows: $\mu'\le \mu$ if and only if the
 convex hull of the Weyl orbit of $\mu$ contains $\mu'$. They then prove:

\begin{main} Let $\mu$, $\mu'$ be points of Atiyah-Bott type for $\xi_0$. If
the closure of a stratum
 $\mathcal{C}_{\mu'}$ meets another stratum $\mathcal{C}_{\mu}$ then
$\mu'\le \mu$ in the Atiyah-Bott ordering.
\end{main}

In particular, it follows that the relation $\preceq$ is a partial ordering
on the set of strata.

The purpose of this paper is to prove a converse to this theorem, by showing:

\begin{main}
Let   $\xi_0$ be a $C^\infty$-principal
$G$-bundle over $C$. Suppose that $\mu,\mu'$ are points of Atiyah-Bott
type for $\xi_0$.
If $\mu'\le \mu$, then $\mathcal{C}_{\mu'}\preceq \mathcal{C}_{\mu}$.
Thus,  $\mathcal{C}_{\mu'}\preceq\mathcal{C}_{\mu}$
if and only if $\mu'\le \mu$ in the
Atiyah-Bott partial ordering.
\end{main}

In general, the decomposition $\{\mathcal{C}_\mu\}_\mu$ of the space of
$(0,1)$-connections is not a
stratification in the sense that the closure of a stratum is not a
union of the higher strata with respect to the Atiyah-Bott ordering. There is
only an inclusion. In Section~\ref{highergenus} we give an example for the
group
$SL(3)$ and any curve of genus greater than one  to show that this inclusion
is not in general an equality. One could ask if a somewhat stronger statement
than the above theorem holds: given $\mu'\le \mu$, does there exist a point
of $\mathcal{C}_\mu$ which is in the closure of $\mathcal{C}_{\mu'}$? It is
possible that the techniques of this paper can be extended to prove this
somewhat stronger statement.

For curves of genus one the situation is much better:
in this case we have a stratification in the strong sense.

\begin{main}
Suppose that $g(C) =1$,
and fix a $C^\infty$ principal $G$-bundle
$\xi_0$ over $C$. Let $\mu$ be a point of  Atiyah-Bott type for $\xi_0$.
Then the closure of the stratum
$\mathcal{C}_\mu$ in the space of $(0,1)$-connections  on $\xi_0$ is
the union $\bigcup_{\mu'\ge \mu}\mathcal{C}_{\mu'}$ where $\mu'$
ranges over points of Atiyah-Bott type for $\xi_0$.
\end{main}

Our motivation for this work came out of the study of the moduli space of
semistable $G$-bundles over an elliptic curve. One can describe this moduli
space as a space of all nontrivial deformations of a ``minimally unstable"
$G$-bundle, which makes clear its structure as a weighted projective space.
For many cases (but not for $SL(n)$), it turns out that there is in fact
essentially a unique such minimally unstable bundle. This paper is our
attempt to understand why this should be so.

The contents of this paper are as follows. In Section 1, we collect some
preliminaries on root systems and parabolic subgroups. In Section 2, we
define points of Atiyah-Bott type and discuss the Harder-Narasimhan
reduction, the strata, and the Atiyah-Bott ordering. Section 3 deals with
what we call harmonic and superharmonic functions on a Dynkin diagram, which
are convenient ways to record some of the positivity properties of the Cartan
matrix. In Section 4, we formulate the main theorems as a combinatorial
problem and show how this problem can be translated into a set of results
about bundles. Most of these results concerning bundles are also established
there. However, one such result requires the notion of an elementary
tranformation of
$G$-bundles. This is defined in Section 5, and we then prove the relevant
bundle result. Finally, in Section 6 we study the problem of finding minimal
elements in the Atiyah-Bott ordering such that the strata correspond to
unstable bundles.

\section{Preliminaries}

\subsection{Basic notation}

We denote by $\frak g$ the Lie algebra of $G$, by $Z(G)$ the center of
$G$ and by
$\frak z_G$ the center of $\frak g$. Fix a maximal torus
$H$ for
$G$ with associated Cartan subalgebra $\frak h$, so that $\frak z_G\subseteq
\frak h$.  There is a
direct sum decomposition $\frak h=(V^*\otimes \Cee)\oplus \frak z_G$. If $K$
is a maximal compact subgroup of $G$, we let $\frak g_\Ar = \sqrt{-1}\Lie
K$, and define $\frak h_\Ar =\frak h\cap \frak g_\Ar$,
$(\frak z_G)_\Ar = (\frak z_G)\cap \frak g_\Ar$. Let
$R$ be the  root system of
$(G,H)$ with Weyl group $W=W(R)$. Fix a set  $\Delta$ of simple roots for
$R$, and let
$R^+$ be the corresponding set of positive roots. The roots
$R$ span a real vector space $V= \operatorname{Ann}(\frak
z_G)_\Ar\subseteq
\frak h^*_\Ar$.  There exists a
$W$-invariant positive definite inner product $\langle \cdot, \cdot
\rangle$ on $V$. Given a root
$\alpha$, there is an associated coroot $\alpha \spcheck\in \frak h_\Ar$.
Using  the inner product to identify $V$ with $V^*$, we have $\alpha
\spcheck = 2\alpha/\langle \alpha, \alpha\rangle$. We denote the
Cartan integer  $ \alpha(\beta \spcheck)$ by $n(\alpha,
\beta)$. Denote by $\Delta\spcheck$ the set of coroots dual to  the simple
roots. The {\sl coroot lattice\/}
$\Lambda$ is the lattice inside $V^*$ spanned by the coroots, and
$\Delta\spcheck$ is a basis for $\Lambda$. The fundamental group
$\pi_1(H)\subseteq \frak h$ and $(2\pi\sqrt{-1})^{-1}\cdot \pi_1(H)$ is a
lattice in $\frak h_\Ar$ containing $\Lambda$. From now on, we will omit the
factor of $(2\pi\sqrt{-1})^{-1}$ and denote the lattice inside $\frak h_\Ar$
as $\pi_1(H)$.  Given
  $\alpha \in \Delta$,  the {\sl fundamental weight\/} $\varpi _{\alpha}\in
V$ is the unique element of ${\frak h}^*$ vanishing on $\frak z_G$ and such
that
$\varpi _{\alpha}(\beta \spcheck) =
\delta_{\alpha\beta}$. The fundamental coweights
$\varpi_{\alpha}\spcheck\in V^*$ are defined similarly. The positive Weyl
chamber $\ov C_0$ in $\frak h_\Ar$ is defined by
$$\ov C_0 = \{ x\in \frak h_\Ar: \alpha (x) \geq 0 \text{ for all $\alpha \in
\Delta$}\}.$$

\subsection{Parabolic subgroups}

The conjugacy classes of parabolic subgroups of $G$ are in natural one-to-one
correspondence with non-empty subsets $I\subseteq \Delta$. The class of
parabolics associated with $I\subseteq \Delta$ has a representative $P^I$ which
is the connected subgroup of $G$ whose Lie algebra consists of the direct sum
of the Cartan subalgebra, all the positive root spaces and the negative root
spaces whose roots are linear combinations of the $\alpha\in
\Delta -  I$. Clearly, $P^I\subseteq P^J$ if and only if $J\subseteq I$.
Thus, the maximal parabolics are of the form $P^{\{\alpha\}}$ for a simple
root $\alpha$. We set $P^{\{\alpha\}} = P^\alpha$. By convention we define
$P^\emptyset=G$.

The Levi factor of $P^I$ is a reductive subgroup $L^I$. Its Lie algebra is
the subalgebra of ${\frak g}$ spanned by ${\frak h}$ and by the root spaces
corresponding to the set of roots in the linear span of
$\Delta -  I$. In particular, $H$ is a maximal torus of $L^I$. Let $\Lambda
_{L^I}\subseteq \Lambda$ be the lattice spanned by the coroots
$\beta\spcheck\in {\frak h}$ dual to the simple roots
$\beta\in\Delta  -  I$. It is a direct summand of $\Lambda$.   The center
$\frak z_{L^I}$ of $\Lie L^I$ is given by
$$\frak z_{L^I} = \bigcap _{\beta \in \Delta - I}\Ker \beta,$$
and hence $\frak z_{L^I}$ is spanned over $\Cee$ by $\varpi_\alpha\spcheck,
\alpha
\in I$ and  $\frak z_G$. Thus $\frak h = \frak z_{L^I} \oplus
(\Lambda_{L^I}\otimes \Cee)$. The following is the equivalent formulation in
terms of root systems:

\begin{lemma}\label{basis} The set $\{\beta\spcheck: \beta \notin I\} \cup
\{\varpi_\alpha\spcheck: \alpha \in I\}$ is a basis for $\Lambda\otimes \Ar
= V^*$.
\end{lemma}
\begin{proof} Suppose that $x\in V$ and that $x$ vanishes on all the
elements in the claim. Since $x$ vanishes on  $\varpi_\alpha\spcheck$ for
$\alpha\in I$, it can be written as a linear combination of the
$\beta\notin I$. Since the
subdiagram of the Dynkin diagram for $G$ spanned by the vertices of
$\Delta-I$ is the Dynkin diagram of a semisimple group, it
follows that $x$ is zero. This shows that
the given elements span $V$. Since the number of them is equal to the
dimension of $V$, it follows that they are a basis.
\end{proof}

A very similar argument shows the following:

\begin{lemma}\label{zL} If $x\in \frak z_L$ is such that $\chi (x) =0$ for
all characters $\chi$ of $G$ and $\varpi_\alpha(x) =0$ for all $\alpha \in
I$, then $x=0$. \endproof
\end{lemma}

A character
$\chi\colon P^I\to \Cee^*$ is {\sl dominant} if its differential
$\chi_*\colon \frak h\to \Cee$ takes nonnegative values on every simple
coroot. The dominant characters of
$P^I$ which vanish on the identity component of the center of $G$
are the nonnegative integral linear combinations of the
fundamental weights $\varpi_\alpha, \alpha \in I$, which  are characters of
$P^I$. For example, if $G$ is simply connected, the dominant characters of
$P^I$ are  exactly the nonnegative integral linear combinations of the
fundamental weights $\varpi_\alpha, \alpha \in I$.

\section{Stratification of the space of all $(0,1)$-connections}

The $C^\infty$ $G$-bundle $\xi_0$ has a first Chern class $c=c_1(\xi_0)\in
H^2(C;\pi_1(G))=\pi_1(G)=\pi_1(H)/\Lambda$, which determines the
topological type of $\xi_0$. Every element of $\pi_1(G)$ arises in this way.
There is a slightly weaker invariant $\zeta(\xi_0)$,   the image
of $c$ in
$\pi_1(H)/\widehat\Lambda$, where
$\widehat\Lambda$ is the saturation of $\Lambda$ in $\pi_1(H)$. The quotient
$\pi_1(H)/\widehat\Lambda$ is a lattice in $\frak h_\Ar/V^*\cong (\frak
z_G)_\Ar$, and we shall view $\zeta(\xi_0)$ as an element of $(\frak
z_G)_\Ar$.

 A holomorphic structure $\xi$ on $\xi_0$ is determined by a
$(0,1)$-connection on $\xi_0$. We denote by $\mathcal{A}=\mathcal{A}(\xi_0)$
the space of all
$(0,1)$-connections on $\xi_0$. It is naturally an affine space for
$\Omega^{(0,1)}(C ;\ad \xi_0)$, and hence supports a natural structure of an
infinite dimensional affine complex variety. Let $\mathcal{G}$ be
the  group of $C^\infty$-automorphisms of $\xi_0$. Different
$(0,1)$-connections determine isomorphic bundles if and only if they differ
by the action of  $\mathcal{G}$ on
$\mathcal{A}$.

A holomorphic $G$-bundle $\xi$ over $C$ is   {\sl semistable} if
$\ad \xi$ is a semistable  vector bundle.
The bundle $\xi$ is semistable if and only if, for every irreducible
representation $\pi \colon G\to GL(N)$, the associated vector bundle $\xi
\times _\pi \Cee^N$ is semistable.  The subset of $\mathcal{A}$
corresponding to semistable bundles is connected, open, and dense  \cite{Ra}.

\subsection{Atiyah-Bott points}

Let $L$ be a reductive subgroup of $G$ containing $H$, and let
$\Lambda_L\subseteq \Lambda$ be the coroot lattice of $L$. Thus $\pi_1(L)=
\pi_1(H)/\Lambda_L$. Let $\frak z_L\subseteq \frak h$ be the Lie algebra of
the center of
$L$. We wish to understand the topological types of
reductions of $\xi_0$ to an $L$-bundle $\eta_0$.  A convenient way to record
this information is through the Atiyah-Bott point of $\eta$.

\begin{defn} Let $L$ be a reductive group and let $H$ be a maximal torus of
$L$. Let $\eta_0$ be a $C^\infty$ $L$-bundle  over
$C$, so that $c_1(\eta_0) \in \pi_1(H)/\Lambda_L$. The {\sl Atiyah-Bott} point
of
$\eta_0$ is the unique point $\mu(\eta_0)\in(\frak z_L)_\Ar$ such that, for
all characters $\chi$ of $L$,   $\chi(\mu(\eta_0))=c_1( \eta_0\times
_\chi\Cee)$. In other words, in the notation at the beginning of this
section, $\mu(\eta_0) = \zeta (\eta_0)$. More generally, suppose that
$P$ is an arbitrary linear algebraic group whose unipotent radical is $U$ and
such that
$P/U=L$, and that $\eta_0$ is a
$P$-bundle. Of course, the bundle $\eta_0$ is topologically equivalent to
$(\eta/U)\times_LP$ for any section of the quotient map $P\to L$.  We define
the Atiyah-Bott point   $\mu(\eta_0)$ to be   $L$-bundle $\mu(\eta_0/U)$.
\end{defn}

\begin{lemma}\label{name}  Suppose that $\eta_0$ is a reduction
of $\xi_0$ to a standard parabolic subgroup
$P^I$ for some $I\subseteq \Delta$, possibly empty. The Atiyah-Bott point
$\mu(\eta_0)$  and the topological type of $\xi_0$ as a
$G$-bundle determine the topological type of $\eta_0/U^I$ as an
$L^I$-bundle {\rm(}and hence of $\eta_0$ as a $P^I$-bundle{\rm)}. Given a
point
$\mu\in
\frak h_\Ar$, there is a reduction of $\xi$ to a $P^I$-bundle whose
Atiyah-Bott point is $\mu$ if and only if the following conditions hold:
\begin{enumerate}
\item[\rm (i)] $\mu$ lies in the Lie algebra $(\frak z_{L^I})_\Ar$ of the
center of
$L^I$.
\item[\rm (ii)] For every simple root $\alpha\in I$ we have
$\varpi_\alpha(\mu) \equiv
\varpi_\alpha(c) \pmod \Zee$.
\item[\rm (iii)] $\chi(\mu)=\chi(c)$ for all characters $\chi$ of $G$.
\end{enumerate}
\end{lemma}

\begin{proof} Let $L=L^I$. Let $\widehat\Lambda_L$ be the saturation of
$\Lambda_L$ in
$\pi_1(H)$, and define $\widehat \Lambda$ similarly. Since $\Lambda_L$ is a
direct summand of $\Lambda$, there is an induced injection
$\widehat\Lambda_L/\Lambda_L\to
\widehat \Lambda/\Lambda$. Suppose that $\alpha\in I$. Then
$\varpi_\alpha\colon
\frak h\to
\Cee$ vanishes on  $\Lambda_L$. Of course, it takes rational values on
$\pi_1(H)$ and integral values on $\Lambda$. Thus,   $\varpi_\alpha$
determines a homomorphism
$\pi_1(H)/\Lambda$ to $\Q/\Zee$, and thus by restriction a homomorphism
$\widehat \Lambda/\Lambda\to\Q/\Zee$. In fact, an easy argument shows that the
sequence
$$0\to \widehat\Lambda_L/\Lambda_L\to \widehat \Lambda/\Lambda\to
\bigoplus_{\alpha\in  I}\Q/\Zee$$
is exact, where the map on the right is the one induced by
$\oplus_{\alpha\in I}\varpi_\alpha$.
The subspace
$\Lambda_L\otimes
\Ar\subseteq
\frak h_\Ar$   is a complementary subspace in $\frak h_\Ar$ to $(\frak
z_L)_\Ar$. Thus projection from
$\frak h_\Ar$ to $(\frak z_L)_\Ar$ defines a homomorphism from
$\pi_1(H)/\Lambda_L$ to $(\frak z_L)_\Ar$, whose kernel is
$\widehat\Lambda_L/\Lambda_L$. Summarizing, we have a commutative diagram
with exact columns and exact first row:
$$\begin{CD}
@. 0 @. 0 @.\\
@. @VVV   @VVV @.\\
0 @>>>\widehat\Lambda_L/\Lambda_L @>>> \widehat \Lambda/\Lambda @>>>
\bigoplus_{\alpha\in  I}\Q/\Zee \\
@. @VVV @VVV @|\\
@. \pi_1(H)/\Lambda_L @>>> \pi_1(H)/\Lambda  @>>>
\bigoplus_{\alpha\in  I}\Q/\Zee\\
@. @VVV @VVV @.\\
@. \pi_1(H)/\widehat\Lambda_L  @>>> \pi_1(H)/\widehat\Lambda @.
\end{CD}$$
The point
$\mu(\eta_0)$ is the image of $c_1(\eta_0)$ under the homomorphism
$\pi_1(H)/\Lambda_L\to (\frak z_L)_\Ar$.
It follows that $\mu(\eta_0)$ determines $c_1(\eta_0)$ up to an element of
$\widehat\Lambda_L/\Lambda$.  Thus, two distinct points $\gamma\not=\gamma'\in
\pi_1(H)/\Lambda_L$ with the same image in $\frak h_\Ar/(\Lambda_L\otimes
\Ar)$ have distinct images in
$\pi_1(H)/\Lambda$. This shows that given $\mu$, there is at most one lift of
$\mu$ to a point $\gamma\in \pi_1(H)/\Lambda_L$ whose image in
$\pi_1(H)/\Lambda$ is
$c$, and hence that
$\mu(\eta_0)$ and $c$ (or equivalently the topological type of $\xi_0$)
determine the topological type of $\eta_0/U$ as an $L$-bundle. Of course, the
topological type of $\eta_0/U$ as an $L$-bundle determines the topological
type of $\eta_0$ as a $P$-bundle.

Next let us show that the Atiyah-Bott point $\mu(\eta_0)$ satisfies the three
conditions stated in the lemma. By construction it lies in the Lie algebra of
the center of $L$.  Clearly, since $c_1(\eta_0)\in
\pi_1(H)/\Lambda_L$ maps to
$c$ in $\pi_1(H)/\Lambda$, the congruence given in the second item holds.
Lastly, since $\chi(\mu(\eta_0))=\chi(c_1(\eta_0))$ for every character
$\chi$ of
$L$, and since under the inclusion of $L\subseteq G$ the element
$c_1(\eta_0)$ maps to $c$, it follows that
$\chi(\mu(\eta_0))=\chi(c)$ for every character of $G$.

Conversely, fix  a point $\mu$ satisfying the three conditions above.  We
shall show that $\mu$ is the Atiyah-Bott point of a reduction to $P^I$ of the
$C^\infty$ bundle
$\xi_0$. We claim that there is a $C^\infty$ $L$-bundle $\eta_0$
such that $c_1(\eta_0)$ projects to $\mu$ and such that
$c_1(\eta_0\times_{L}G)=c$.
It suffices to show that there is an element $c_{L}\in
\pi_1(H)/\Lambda_{L}$ which projects to $\mu$ and whose image in
$\pi_1(H)/\Lambda$ is
$c$, for then there is a $C^\infty$ $L$-bundle $\eta_0$ with
$c_1(\eta)= c_L$, and this bundle has the required properties. Choose
$\widetilde c\in \pi_1(H)$ lifting $c$. For all $\alpha \in I$,
$\varpi_\alpha (\mu - \widetilde c) \in \Zee$. Thus there exists a $\lambda
\in \Lambda$ such that $\varpi_\alpha (\mu) = \varpi_\alpha(\lambda +
\widetilde c)$ for all $\alpha \in I$. Let $c_L$ be the image of $\lambda +
\widetilde c$ in $\pi_1(H)/\Lambda_{L}$. Then $c_L$ maps to $c\in
\pi_1(H)/\Lambda$, and the image $\mu'$ of $c_L$ satisfies:
$\varpi_\alpha(\mu) = \varpi_\alpha(\mu')$ for all $\alpha\in I$. Clearly,
for every character $\chi$ of $G$, $\chi(\mu) = \chi(\mu')$. By
Lemma~\ref{zL}, $\mu =\mu'$, and so $c_L$ is the required element.  This
completes the proof of Lemma~\ref{name}.
\end{proof}

\begin{defn}  A pair $(\mu,I)$
consisting of a point
$\mu\in \frak h_\Ar$ and a subset $I\subseteq \Delta$ is said to be {\sl of
  Atiyah-Bott type for
$c$} (or for $\xi_0$) if and only if the three conditions given in
Lemma~\ref{name} are satisfied. If
$(\mu,I)$ and $(\mu,I')$ are Atiyah-Bott pairs of type $c$, then so is
$(\mu,I\cap I')$. Thus, one can always choose
$I$ to be minimal so that the conditions hold. This minimal $I$ consists of
all $\alpha\in\Delta$ for which $\alpha(\mu)\not= 0$. A point $\mu\in
\frak h_\Ar$ is said to be {\sl of Atiyah-Bott type for $c$} if there is a
subset
$I\subseteq \Delta$ such that $(\mu,I)$ is a pair of Atiyah-Bott type for $c$.
\end{defn}

\begin{corollary}  Given a pair
$(\mu,I)$ of Atiyah-Bott type for $c$, there is a reduction of $\xi_0$
to an $L^I$-bundle $\eta_0$ whose Atiyah-Bott point is $\mu$, and $\eta_0$ is
unique up to $C^\infty$ isomorphism.
\endproof
\end{corollary}

\subsection{The Harder-Narasimhan reduction}

Recall that every unstable holomorphic $G$-bundle $\xi$ has a
 canonical reduction of its structure group
to a conjugacy class of parabolic subgroups. We can always choose a
standard parabolic $P^I$ in this class. We call this reduction the
{\sl Harder-Narasimhan reduction\/} and $P$ the {\sl Harder-Narasimhan
parabolic\/} of $\xi$. If $\xi$ is semistable, then by convention we take
as the Harder-Narasimhan reduction the trivial reduction to the group $G =
P^\emptyset$. The following summarizes the basic properties of this
reduction~\cite{Ra},
\cite{AtBo}, \cite[Corollary 2.11]{FM}.

\begin{lemma}\label{name1} Let $\xi$ be a holomorphic $G$-bundle over
$C$. A reduction of $ \xi$ to a $P^I$-bundle  $\eta$ is the
Harder-Narasimhan reduction if and only if
\begin{enumerate}
\item[\rm (i)] $\eta/U^I$ is a semistable  $L^I$-bundle.
\item[\rm (ii)] $\mu(\eta)\in \ov C_0$.
\item[\rm (iii)] If $\alpha\in I$, then
$\alpha(\mu(\eta))>0$. \endproof
\end{enumerate}
\end{lemma}

\begin{corollary} Let
$\mu\in \frak h_\Ar$ be of Atiyah-Bott type for $c$ and
let $I=\{\alpha \in \Delta: \alpha(\mu)>0\}$.  Then there is a
holomorphic $G$-bundle structure on $\xi_0$ which has Harder-Narasimhan
reduction to a subgroup $P^I$ with Atiyah-Bott point $\mu$ if and only if
$\mu \in \ov C_0$.
\end{corollary}

\begin{proof} Since every $C^\infty$ bundle has a semistable  holomorphic
structure, this is clear from the previous lemma and
Lemma~\ref{name}.
\end{proof}

\subsection{The strata}

Following Shatz and Atiyah-Bott we define a stratification of the space of
$(0,1)$-connections on $\xi_0$.

\begin{defn} Fix a  point $\mu\in\ov C_0$ of Atiyah-Bott
type for $c$. The {\sl stratum\/}
$\mathcal{C}_\mu\subseteq \mathcal{A}$ is
the set of all $(0,1)$-connections defining holomorphic structures on
$\xi_0$ whose Harder-Narasimhan reduction has Atiyah-Bott point equal to
$\mu$.
\end{defn}

In particular, if $c=\zeta+v$ with $\zeta\in\frak z_G$ and $v\in V^*$, then
$\zeta$ is of Atiyah-Bott type for $c$ and
$\mathcal{C}_\zeta$ is exactly the open dense set of $(0,1)$-connections
defining semistable  $G$-bundles. Thus for example if $G$ is semisimple,
then the stratum of semistable bundles is $\mathcal{C}_{0}$.

The
strata are preserved by the action of  $\mathcal{G}$.
 The union of these strata over all $\mu\in \ov C_0$ of Atiyah-Bott type for
$\xi_0$ is $\mathcal{A}$.

 We say that a holomorphic bundle structure on
$\xi_0$, or equivalently a holomorphic bundle whose underlying topological
bundle is isomorphic to $\xi_0$, is contained in the stratum
$\mathcal{C}_\mu$ if all
$(0,1)$-connections determining holomorphic bundles
isomorphic to it  are in this stratum.

\begin{proposition}
Suppose that $\mu\in \ov C_0$
is a point of Atiyah-Bott type for $c$. Let $a,b\in \mathcal{C}_\mu$. Then
there is a connected complex space $S$, a holomorphic family of
$(0,1)$-connections $A_s, s\in S$, and points $s_1, s_2\in S$ such that
$A_s\in \mathcal{C}_\mu$ for all $s\in S$ and
$A_{s_1}$ is gauge equivalent to $a$ and $A_{s_2}$ is gauge equivalent to
$b$. In particular, $\mathcal{C}_\mu/\mathcal{G}$ is
connected.
\end{proposition}

\begin{proof}
 Lemma~\ref{name} shows that the parabolic subgroup and the
topological type of the $L$-bundle $\eta_0$ are determined by
$\mu$ and the topological type of $\xi_0$.
Given two connections $a_L$, $b_L$ on $\eta_0$ which define semistable
holomorphic structures, an open dense set of the line in
$\mathcal{A}(\eta_0)$ joining them will also define semistable holomorphic
structures. We can use this line to join the corresponding
$(0,1)$-connections $a'$ and $b'$ on $\xi_0$. The space of holomorphic
$P$-bundles with a
given $L$-reduction is an affine space,
by~\cite[Appendix]{FMII}. After choosing a $C^\infty$ trivialization of these
spaces, we can find holomorphic, connected families of $(0,1)$-connections
joining $a'$ to $a$ and $b'$ to $b$.
\end{proof}

The strata are locally closed in the Zariski
topology in the following
sense. If $T$ is a finite dimensional parameter space for an algebraic
family of holomorphic structures on $\xi$ and $\mu$ is a point
of Atiyah-Bott type for $\xi$, then the subspace of $T$ consisting of
points parametrizing bundles contained in $\mathcal{C}_\mu$ is a
locally closed subspace of $T$ with respect to the Zariski topology.
An analogous statement for the classical topology holds for analytic
families of holomorphic structures on $\xi$.

\begin{defn}
Following
 Atiyah-Bott,
we say that
$\mathcal{C}_{\mu_1} \preceq _r \mathcal{C}_{\mu_2}$ if there exists a
holomorphic $G$-bundle $\xi$ in $\mathcal{C}_{\mu_2}$ and an arbitrarily
small deformation of it to a bundle in $\mathcal{C}_{\mu_1}$. In light
of the above remarks this is equivalent to the existence of a
holomorphic family of
$G$-bundles $\Xi$ over $C\times T$, where $T$ is connected, and a point
$t_0\in T$, such that $\Xi _t\in \mathcal{C}_{\mu_1}$ for $t\not= t_0$, and
such that $\Xi _{t_0} \in \mathcal{C}_{\mu_2}$. The relation
$\preceq_r$ generates a transitive relation  which we
denote by $\preceq$.
\end{defn}

\subsection{The Atiyah-Bott ordering}

\begin{defn}\label{ABorder} We define the {\sl Atiyah-Bott partial ordering\/}
on
  $\frak h_\Ar$ as follows. Given $x\in \frak h_\Ar$, let $\widehat{W\cdot x}$
denote the convex hull of the finite set
$W\cdot x$. For
$x, y\in \frak h_\Ar$, we define $x\geq y$ if and only if $y \in
\widehat{W\cdot
  x}$.  Notice that $x\geq y$ if and only if the projections of $x$ and
$y$ into $(\frak z_G)_\Ar$ are equal and
$v^*(x)\geq v^*(y)$ where $v^*(x)$ and $v^*(y)$ are the projections of $x$ and
$y$ into $V^*$.
\end{defn}

For points in $\ov C_0$, there is a simple
characterization of the ordering.

\begin{lemma}\label{orderchar} Suppose that $x,y\in \ov C_0$. Then
  $x\geq y$ if and only if:
\begin{enumerate}
\item[\rm (i)]  For every simple root $\alpha$, $\varpi _\alpha(x) \geq
\varpi_\alpha (y)$.
\item[\rm (ii)] The projections of $x$ and $y$ into $\frak z_G$ are equal.
\end{enumerate}
\end{lemma}

\begin{proof} Since the Weyl group acts trivially on $\frak z_G$, it suffices
to divide out by
  $\frak z_G$. We may thus assume that $G$ is semisimple. By
\cite[Lemma 12.14]{AtBo}, for $x,y\in \ov C_0$, $x\geq
  y$ if and only if $\langle t, x\rangle \geq \langle t, y\rangle$ for
  all $t\in \ov C_0$. But the simplicial cone $\ov C_0$ is spanned over
  $\Ar^+$ by the elements $\varpi_\alpha\spcheck$, $\alpha \in
  \Delta$, and $\langle \varpi_\alpha\spcheck, x\rangle =
  c\varpi_\alpha(x)$ for some positive constant $c$. Thus $x\geq y$ if
  and only if $\varpi _\alpha(x) \geq \varpi_\alpha (y)$  for every
  $\alpha\in \Delta$.
\end{proof}

The next corollary says that on $\ov C_0$
the Atiyah-Bott partial ordering really is a partial ordering:

\begin{corollary} Suppose that $x_1, x_2\in \ov C_0$ are such that $x_1\leq
x_2$ and $x_2 \leq x_1$. Then $x_1=x_2$. More generally, given $x_1, x_2\in
\frak h_\Ar$, $x_1\leq x_2\leq x_1$ if and ony if $x_1$ and $x_2$ are
conjugate under $W$. \endproof
\end{corollary}

The relevance of the Atiyah-Bott ordering to the ordering of strata in the
space of $(0,1)$-connections is given by the following theorem
\cite{AtBo}:

\begin{theorem}\label{mainthm1} {\bf (Atiyah-Bott)}
Suppose that $\mathcal{C}_{\mu_1} \preceq_r
  \mathcal{C}_{\mu_2}$, i.e., suppose that there is a bundle in the
  stratum $\mathcal{C}_{\mu_2}$ and an arbitrarily small deformation
  of it which is a bundle in the stratum
  $\mathcal{C}_{\mu_1}$. Then $\mu_1\leq \mu_2$. \endproof
\end{theorem}

\begin{corollary}
If $\mathcal{C}_{\mu_1}\preceq \mathcal{C}_{\mu_2}$, then
$\mu_1\leq\mu_2$ in the Atiyah-Bott partial ordering.
In particular, we see that the transitive
relation $\preceq$ is a partial ordering on
the set of strata $\mathcal{C}_\mu$.
\end{corollary}

The rest of this paper is devoted to establishing a converse to this
result. For the case of curves of genus one, we have the following strong
converse:

\begin{theorem}\label{mainthm}  Suppose that $g(C) =1$, and that $\mu_1\leq
\mu_2$ in the Atiyah-Bott ordering and that $\xi$ is a holomorphic bundle in
$\mathcal{C}_{\mu_2}$. Then there is an arbitrarily small deformation $\xi'$
of $\xi$  contained in $\mathcal{C}_{\mu_1}$.
\end{theorem}

For higher genus, there is a weaker version of Theorem~\ref{mainthm}:

\begin{theorem}\label{mainthm2} Let $C$ be a smooth curve of genus
at least one. Then $\mathcal{C}_{\mu_1}\preceq
  \mathcal{C}_{\mu_2}$, if and only if $\mu_1\leq\mu_2$.
\end{theorem}

It follows from Theorem~\ref{mainthm} that
in the case of curves of genus one,
the stratification is in fact a
stratification in a strong sense as the next corollary shows.

\begin{corollary} Suppose that $g(C) =1$.   Let $\mu\in \ov C_0$ be a point
of Atiyah-Bott type for $c$. The closure of the stratum $\mathcal{C}_\mu$ in
$\mathcal{A}$ is equal to $\bigcup_{\mu'\ge
  \mu}\mathcal{C}_{\mu'}$ where $\mu'$ ranges over points in $\ov C_0$ of
Atiyah-Bott type for $c$.
\end{corollary}

\begin{proof}   Suppose that $\mu'\in \ov C_0$ is a point of
Atiyah-Bott type  with
$\mu'\ge \mu$. Let $A$ be a $(0,1)$-connection in $\mathcal{C}_{\mu'}$. Let
$\eta'$ be the holomorphic $G$-bundle determined by $A$. According to
Theorem~\ref{mainthm} there is an arbitrarily small deformation of $\eta'$ to
a holomorphic $G$-bundle
$\eta$ contained in $\mathcal{C}_\mu$. We can extend $A$ to a
$(0,1)$-connection on this deformation. A $C^\infty$ trivialization of the
deformation allows us to view all the $(0,1)$-connections in the family as
connections on   $\xi_0$. The resulting $(0,1)$-connection determining
$\eta$ is then arbitrarily close to
$A$ in the space of $(0,1)$-connections. Thus, every neighborhood
of $A$ in $\mathcal{A}$ contains a $(0,1)$-connection in $\mathcal{C}_{\mu}$.
Hence, $A$ is contained in the closure of $\mathcal{C}_\mu$.

This proves that the subspace $\bigcup_{\mu'\ge
  \mu}\mathcal{C}_{\mu'}$, where $\mu'$ ranges over points in $\ov C_0$ of
Atiyah-Bott type for $c$, is contained in the closure of $\mathcal{C}_\mu$.
The Atiyah-Bott result is that the closure of $\mathcal{C}_\mu$ is contained
in this union. Hence, the closure of $\mathcal{C}_\mu$ is equal to this union.
\end{proof}

\section{Harmonic  functions on a Dynkin diagram}

The Cartan matrix $\bigl(n(\alpha, \beta)\bigr)_{\alpha,\beta\in\Delta}$  has
two fundamental properties: Its off-diagonal entries are non-positive, and its
inverse is a positive matrix.
The fact that these hold for $\Delta$ and for all non-empty subsets of
$\Delta$ allows us to establish a theory of harmonic functions on
$\Delta\spcheck$ with results paralleling those for harmonic  functions on a
compact manifold.

\subsection{The basic definitions}

A point $x\in V^*$ can be written uniquely
$$x=\sum_{\alpha\in \Delta}r_\alpha\alpha\spcheck$$ with $r_\alpha\in\Ar$.
Thus, we can view $x$ as a function $f_x$ on $\Delta\spcheck$ by
$f_x(\alpha\spcheck)=r_\alpha=\varpi_\alpha(x)$. This determines a linear
isomorphism between
$V^*$ and the space of functions
$\Delta\spcheck\to\Ar$.

\begin{defn}\label{order} Let  $f\colon \Delta\spcheck\to \Ar$ be a function.
It is {\sl harmonic at $\alpha\spcheck\in\Delta\spcheck$} if
$$\sum_{\beta\in\Delta}n(\alpha,\beta)f(\beta\spcheck)=0.$$ The
function $f$ is   {\sl superharmonic at
$\alpha\spcheck$} if
$$\sum_{\beta\in\Delta}n(\alpha,\beta)f(\beta\spcheck)\ge 0,$$  and
{\sl subharmonic at
$\alpha\spcheck$} if
$$\sum_{\beta\in\Delta}n(\alpha,\beta)f(\beta\spcheck)\le 0.$$
Thus, for $x\in V^*$, $f_x$ is superharmonic, harmonic, or subharmonic
at $\alpha$ according  to whether  $\alpha(x)\ge 0$, $\alpha(x)=0$, or
$\alpha(x)\le 0$.
Let $A\subseteq \Delta\spcheck$. Then
$f$ is  {\sl harmonic except at $A$} if it is harmonic at every
$\alpha\spcheck\in \Delta\spcheck -  A$. It is   {\sl
harmonic} if it is harmonic at every
$\alpha\spcheck\in \Delta\spcheck$. Similarly, one defines the notions of
superharmonic and subharmonic except at $A$ and superharmonic and subharmonic.
\end{defn}

There is a more geometric way of viewing the harmonic condition at $\alpha\in
\Delta\spcheck$. It is a local condition just involving the values of $f$ on
the unit disk centered at $\alpha\spcheck$ -- an average value condition
where the  values  are weighted by the Cartan integers.
 For $\alpha\spcheck \in \Delta\spcheck$, let $\mathbf{s}(\alpha\spcheck)
  =\{\beta\spcheck\in \Delta\spcheck: n( \alpha, \beta) < 0\}$ be the
``sphere of radius one" around
$\alpha\spcheck$. Since the only positive Cartan integers are the diagonal
ones which are equal to $2$, for a function $f\colon \Delta\spcheck \to \Ar$
it is immediate from the definition that $f$ is superharmonic at
$\alpha\spcheck$ if and only if
$$ 2f(\alpha\spcheck)\ge -\sum _{\beta\spcheck \in
\mathbf{s}(\alpha)}n(\alpha, \beta)f(\beta\spcheck). $$
There are similar descriptions of when $f$ is harmonic or subharmonic at
$\alpha\spcheck$.

\subsection{Basic properties of harmonic functions on $\Delta\spcheck$}

\begin{defn} If $f$ and $g$ are  functions from $\Delta\spcheck$ to $\Ar$ we
say that $f\ge g$ if $f(\alpha\spcheck)\ge g(\alpha\spcheck)$ for all
$\alpha\spcheck\in \Delta\spcheck$.
\end{defn}

\begin{lemma}\label{inC0} For $x\in V^*$, the corresponding function $f_x$ is
superharmonic if and only if $x\in\ov C_0$.
\end{lemma}

\begin{proof} As we have seen, $f_x$ is superharmonic at $\alpha\spcheck$ if
and only if $\alpha (x) \geq 0$. Consequently, $f_x$ is superharmonic if and
only if $x\in \ov C_0$.
\end{proof}

\begin{lemma}\label{superpos} A function $f\colon \Delta\spcheck\to \Ar$ is
superharmonic if and only if it is a nonnegative linear combination of
$\{f_{\varpi_\alpha\spcheck}\}_{\alpha\in \Delta}$. A superharmonic
function $f$ on $\Delta\spcheck$ satisfies
$f\ge 0$. If $f$ is harmonic, then it is zero.
\end{lemma}

\begin{proof} Since the simplicial cone $\ov C_0$ is spanned by
$\{\varpi_\alpha\spcheck\}_{\alpha\in \Delta}$, the first statement
is clear. Since $\varpi_\beta(\varpi_\alpha\spcheck) \geq 0$ for all $\alpha,
\beta \in \Delta$, by e.g.~\cite[p.\ 168]{Bour}, it follows that, if $f$
is superharmonic, then
$f\ge 0$. Finally, if $f$ is harmonic, then $f$ and $-f$ are both
superharmonic, so that $f=0$.
\end{proof}

\begin{lemma}\label{inequal} Let $A\subseteq \Delta\spcheck$, and let $f$,
$g$, and $h$ be functions from $\Delta\spcheck$ to
$\Ar$ with $f$ harmonic except at $A$, $g$ superharmonic except at $A$, and
$h$ subharmonic except at $A$. Suppose that $g(\alpha\spcheck)\ge
f(\alpha\spcheck)\ge h(\alpha\spcheck)$ for all
$\alpha\spcheck\in A$. Then
$g\ge f\ge h$.
\end{lemma}

\begin{proof} Since the negative of a harmonic function is harmonic and the
negative of a superharmonic function is subharmonic, the result for a
subharmonic function $h$ with $f\ge h$ follows from that for a
superharmonic function $g$ with $g\ge f$. After replacing $g$ by $g-f$, we may
assume that
$g$ is superharmonic except at $A$ and that $g|A\ge 0$, and wish to show that
$g\ge 0$. Suppose that   $g=f_x$ for $x\in V^*$ with $x =\sum _{\beta
\spcheck\in \Delta\spcheck}x_\beta\beta\spcheck$, where the $x_\beta$ are
given nonnegative real numbers for $\beta\spcheck \in A$. The function $g$ is
superharmonic except at $A$ if and only if
$\alpha(x) \geq 0$ for all $\alpha\spcheck \notin A$ if and only
$$\sum _{\gamma\spcheck \in \Delta\spcheck -  A }x_\gamma n(\alpha, \gamma)
\geq - \sum _{\beta\spcheck \in A}x_\beta n(\alpha, \beta)=v_\alpha,$$
say, where the $v_\alpha$ are given nonnegative real numbers since $n(\alpha,
\beta)\leq 0$. It suffices to show that $x_\gamma \geq 0$ for all
$\gamma\spcheck \notin A$. This follows since the inverse of the
Cartan matrix for the subdiagram corresponding to
$\Delta\spcheck -  A$  has nonnegative entries.
\end{proof}

\begin{lemma}\label{exten} Fix $A\subseteq \Delta$.
Given  a
function $f_0\colon A\to\Ar$,  there is a unique extension  $f\colon
\Delta\spcheck\to \Ar$ of $f_0$ which is harmonic except at $A$. If
$f_0\ge 0$, then $f\ge 0$.
\end{lemma}

\begin{proof} By
Lemma~\ref{basis} appplied to the dual root system, the set
$\{\varpi_\beta\}_{\beta\in  A}\cup
(\Delta  -  A )$ is a  basis for $V$. Thus
if $f=f_x$ for $x\in V^*$, then $x$ is uniquely determined by the
conditions
$\alpha(x) = 0$ for $\alpha \notin A$ and $\varpi_\beta(x) =
f(\beta\spcheck)$ for $\beta \in A$.  The positivity statement is the special
case of the previous lemma, applied to the inequality $f\geq 0$ on $A$ and
viewing $f$ as superharmonic except at $A$ and $0$ as harmonic.
\end{proof}

\begin{lemma}\label{compare}
If $\mu,\mu'\in\ov C_0$, then
$\mu\ge \mu'$ in the Atiyah-Bott ordering if and only if $f_\mu\ge
f_{\mu'}$. If $\mu,\mu'\in -\ov C_0$, then $\mu\ge \mu'$ in the
Atiyah-Bott ordering if and only if $f_{\mu'}\ge f_{\mu}$.
\end{lemma}

\begin{proof}
The first statement is immediate from Lemma~\ref{orderchar}.
Suppose $\mu,\mu'\in-\ov C_0$. Then $\mu\ge\mu'$ in the Atiyah-Bott
ordering if and only if $\mu'\in \widehat {W\cdot \mu}$, if and only if
$-\mu'\in \widehat{ W\cdot(- \mu)}$. Since $-\mu,-\mu'\in\ov C_0$, it
follows from the first statement that this holds if and only if
$f_{-\mu}\ge f_{-\mu'}$ which is clearly equivalent to the statement
$f_\mu\le f_{\mu'}$.
\end{proof}

\subsection{Examples}

\begin{figure}
\centerline{\scalebox{.70}{\includegraphics{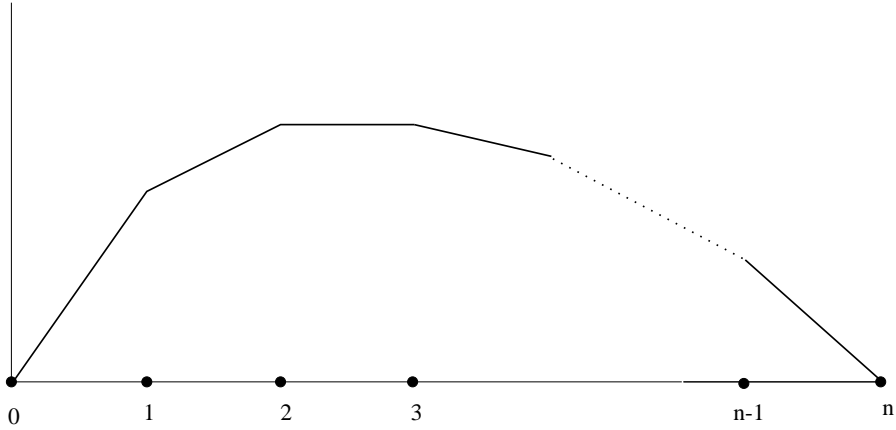}}}
\caption{The Graph of a Superharmonic Function for $A_{n-1}$}\label{fig1}
\end{figure}

Let us suppose that $G=SL(n)$. Then we
view $\Delta\spcheck$ as the points $\{1,\ldots,n-1\}$ in the
interval $I=[0,n]$ with
\begin{equation}\label{dist}
|a-b|=1 {\rm \ \  if\  and\  only\  if\  \ }n(a,b)=-1 {\rm \ \ for\
 all\  } a, b\in \Delta\spcheck.
\end{equation}
Given
a function $f\colon \Delta\spcheck\to\Ar$ we
extend it to a continuous
function $\hat f\colon I\to \Ar$ by first requiring
that $\hat f(0)=\hat f(n)=0$ and that $\hat f$ is linear on each
inteval $[k,k+1]$, $k=0,\ldots,n-1$. Then $f$ is superharmonic if and
only if $\hat f$ is a convex function on $I$. Furthermore, $f$ is harmonic
at $a\in\Delta\spcheck$ if and only if $\hat f$ is linear at $a$.
 (See Figure 1.)

As Atiyah-Bott point out, $f\ge g$ in their ordering if and only if
the graph of $\hat f$ lies above that of $\hat g$.

\begin{figure}
\centerline{\scalebox{.65}{\includegraphics{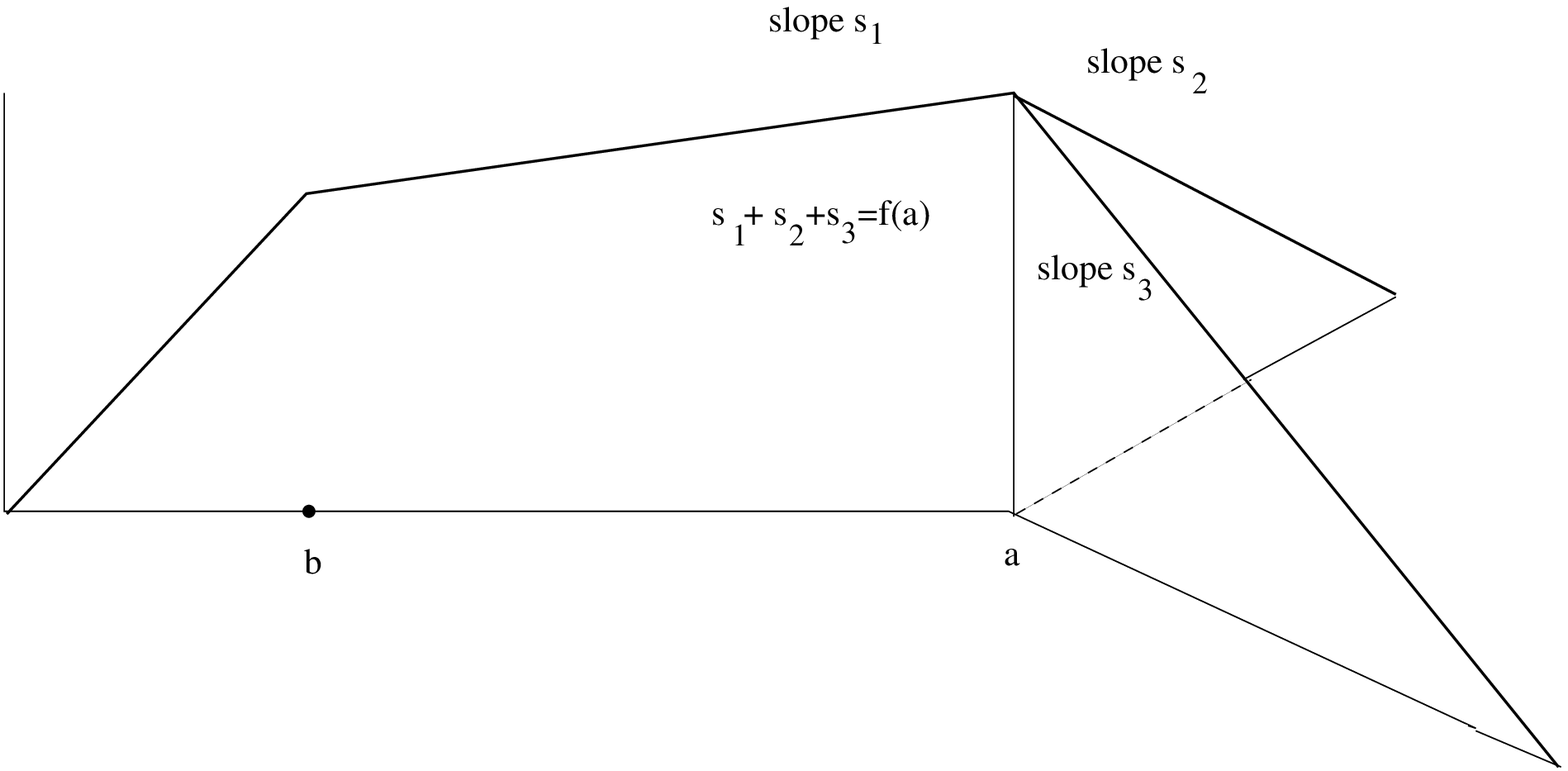}}}
\caption{The Graph of a Function Superharmonic at $b$ and Harmonic Elsewhere}
\end{figure}

Now suppose that $G$ is a simple group of type $D_n$ or $E_6,E_7,$ or
$E_8$, so that the Dynkin diagram for $\Delta\spcheck$ is a union of
three simple
chains (Dynkin diagrams of $A$-type) $\Delta\spcheck_1,\Delta\spcheck_2,
\Delta\spcheck_3$ meeting at the trivalent vertex $a$ of $\Delta\spcheck$.
For $i=1,2,3$ let $\ell_i=\#\Delta\spcheck_i$ and let $I_i$ be the
interval $[0,\ell_i]$. Identify $\Delta\spcheck_i$ with $\{1,2,\ldots,\ell_i\}
\subset I_i$ satisfying Condition~\ref{dist} in such a way that
the trivalent vertex $a$ is identified with $\ell_i\in I_i$. Let
$T=\bigcup_{i=1}^3
I_i$ where $I_i\cap I_j=\{a\}$ for
all $i\not= j$. There is a unique embedding of $\Delta\spcheck\subset
T$ consistent with the given embeddings of $\Delta\spcheck_i\subset I_i$.
Given a function $f\colon\Delta\spcheck\to\Ar$ we extend it to a
continuous function
$\hat f\colon T\to\Ar$ by requiring that $\hat f$ vanishes on $0\in I_i$,
for $i=1,2,3$ and by requiring that $\hat f$ is linear on each interval
of the form $[k,k+1]\subset I_i$, $k=0,1,\ldots,\ell_i-1$. Let $s_i = \hat
f(a) -
\hat f(\ell_i-1)$ be the ``slope" at $a$. The function $f$ is superharmonic if
and only if (i)
$\hat f|I_i$ is convex for each $i=1,2,3$, and
(ii) $\sum_{i=1}^3s_i\ge \hat f(a)$.
The function $f$ is harmonic at $b\in\Delta\spcheck_i-\{a\}$ if
and only if $\hat f|I_i$ is linear at $b$. The function $f$ is
harmonic at $a$ if and only if the inequality in (ii) above is
an equality.
 (See Figure 2.)

\begin{figure}
\centerline{\scalebox{.65}{\includegraphics{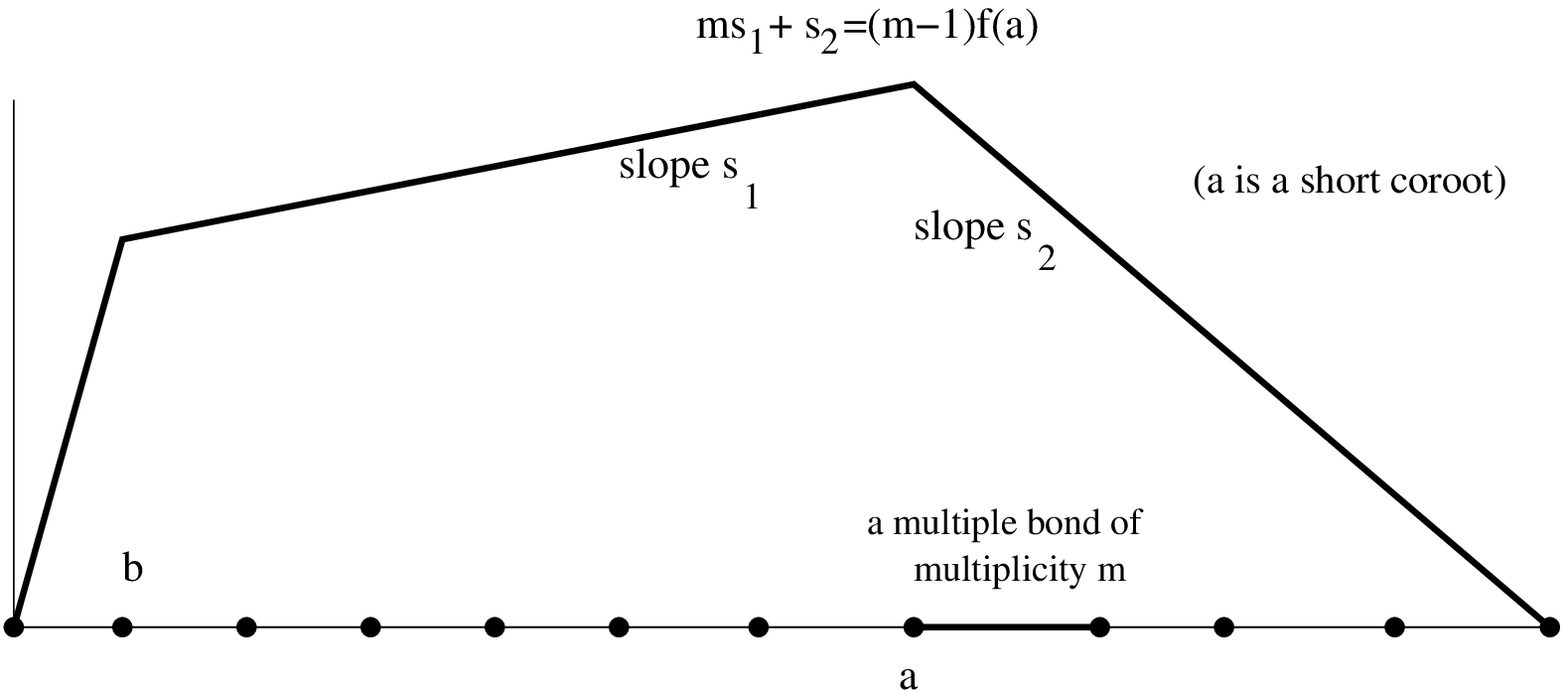}}}
\caption{The Graph of a Function Superharmonic at $b$ and Harmonic Elsewhere}
\end{figure}

Lastly, suppose that $G$ is a simple group of type $B_n,C_n, F_4$ or $G_2$
so that the Dynkin diagram of $G$ is a chain with a single multiple bond.
Identify $\Delta\spcheck$ with the points $\{1,\ldots, n\}$ in the inteval
$I=[0,n+1]$ satisfying Condition~\ref{dist}. Suppose that the multiple
bond in the Dynkin diagram connects
$\ell$ and $\ell+1$ with $\ell$ identified with a short coroot $a$ (which
then is then the special vertex)
and $\ell+1$ being identified with a long coroot
in $\Delta\spcheck$.
Let $m\ge 2$ be the multiplicity of the bond.
Let $f\colon\Delta\spcheck\to\Ar$ be a function. We extend $f$ to a function
$\hat f\colon I\to\Ar$ by requiring that $\hat f(0)=\hat f(n+1)=0$
and that $\hat f$ be linear on each inteval of the form $[k,k+1]$,
$k=0,\ldots,n$. Let $s_1 = \hat f(\ell) - \hat f(\ell -1)$ and let $s_2 =
\hat f(\ell) - \hat f(\ell +1)$. Then
$f$ is superharmonic if and only if (i)
$\hat f|[0,\ell]$ and $\hat f|[\ell,n+1]$ are convex and (ii)
$(m-1)f(a)\le ms_1+s_2$.
The function $f$ is harmonic at $b\in\Delta\spcheck-\{a\}$ if and only if
$\hat f$ is linear at $b$, and $f$ is harmonic at $a$ if and only if
the inequality in (ii) above is an equality.
 (See Figure 3.)

\subsection{Relation to points of Atiyah-Bott type}

We now link the above discussion with points of Atiyah-Bott type.  Recall
that the strata
$\mathcal{C}_\mu$ of $\mathcal{A}=\mathcal{A}(\xi)$ are indexed by
points $\mu$ of Atiyah-Bott type for $c$ and such that $\mu\in \ov
C_0$.

\begin{proposition}\label{abpoint} Write
$c=\zeta+v$ with $\zeta\in(\frak z_G)_\Ar$ and $v\in V^*$. Let  $\mu\in \frak
h_\Ar$. Decompose $\mu=\zeta'+\nu$, with
$\zeta'\in(\frak z_G)_\Ar$ and $\nu\in V^*$.  Let $I\subseteq \Delta$. Then
$(\mu,I)$ is of Atiyah-Bott type for
$c$ if and only if
\begin{enumerate}
 \item[\rm (i)] The function $f_\nu$ is harmonic except at $I$;
\item[\rm (ii)] For each $\alpha\spcheck\in I$,
$f_\nu(\alpha\spcheck)\equiv \varpi_\alpha(c) \pmod \Zee$.
\item[\rm (iii)] $\zeta'=\zeta$.
\end{enumerate}
Moreover, $\mu\in \ov C_0$ if and only if $f_\nu$ is superharmonic.
\end{proposition}

\begin{proof} The condition $\zeta'=\zeta$ is equivalent to the condition
that for each character $\chi$ of $G$ we have $\chi(\mu)=\chi(c)$. The
condition that $f_\nu$ is harmonic at $\alpha\spcheck$ means that
$\alpha(\mu)=0$. The condition that $f_\nu(\alpha\spcheck)\equiv
\varpi_\alpha(c)
\pmod \Zee$ is equivalent to $\varpi_\alpha(\mu )\equiv
\varpi_\alpha(c) \pmod \Zee$. Thus, by Lemma~\ref{name},
$(\mu,I)$ is of Atiyah-Bott type for $c$. The final statement follows from
Lemma~\ref{inC0}.
\end{proof}

\section{Proof of the main theorems}

In this section we prove the two main theorems of this paper.
 In the first subsection, we describe a combinatorial
problem which involves moving from one Atiyah-Bott point of type $c$ in $\ov
C_0$ to a smaller one. In the second, we relate this combinatorial problem to
operations on bundles. Taken together, these results give a proof of
Theorem~\ref{mainthm}. We then describe the modifications necessary to prove
the weaker result Theorem~\ref{mainthm2} in the case of higher genus.

\subsection{A combinatorial discussion}

Let $\mu\in \frak h_\Ar$.
Then we define the function $f_\mu\colon\Delta\spcheck\to \Ar$
as follows. We decompose $\mu=\zeta'+\nu$ with $\zeta'\in(\frak z_G)_\Ar$
and $\nu\in V^*$ and we set $f_\mu=f_\nu$.
Note that  $\mu\in \ov C_0$ if and only if $f_\mu$ is superharmonic.
For points $\mu,\mu'\in\ov C_0$ we have $f_\mu\ge f_{\mu'}$ if and
only if $\mu\ge \mu'$
in the ordering given in Definition~\ref{ABorder}.
Given $f_\mu \geq f_{\mu'}$, the following
theorem gives a description of a sequence of moves which begins with $\mu$
and ends with $\mu'$. Although we shall primarily be interested in the case
where $\mu $ and $\mu'$ lie in $\ov C_0$, we cannot assume that the
intervening points $\mu_i$ lie in $\ov C_0$ and shall not make any
assumptions on $\mu $ and $\mu'$.

\begin{theorem}\label{decrease}  Suppose that we have two
pairs $(\mu,I)$ and $(\mu',I')$
 of Atiyah-Bott  type for $c$ with
$f_\mu\ge f_{\mu'}$. Then there
is a sequence of pairs of Atiyah-Bott  type for $c$,
$(\mu, I)=(\mu_0,I_0), \dots ,(\mu_n, I_n)=(\mu',I')$ with $f_{\mu_i}\ge
f_{\mu_{i+1}}$ for $0\le i\le n-1$, and such that, for each $i,\ 0\le
 i\le n-1$, $(\mu_{i+1},I_{i+1})$ is obtained from $(\mu_i,I_i)$ by
 one of the following moves.
\begin{itemize}
\item[\rm (1)] $I_{i+1}=I_i\cup \{a\}$,
$f_{\mu_i}|_{I_i}=f_{\mu_{i+1}}|_{I_i}$, and $f_{\mu_{i+1}}$ is subharmonic at
$a$. Furthermore, if $\nu\in \frak h_\Ar$ is a point of Atiyah-Bott type
for $c$ and
$f_{\mu_i}\ge f_\nu\ge f_{\mu_{i+1}}$ then  either $f_\nu=f_{\mu_i}$ or
$f_\nu=f_{\mu_{i+1}}$.
\item[\rm (2)] $I_{i+1}\subseteq I_i$,
  $f_{\mu_i}|_{I_{i+1}}=f_{\mu_{i+1}}|_{I_{i+1}}$.
\item[\rm (3)] $I_i=I_{i+1}$ and there is $a\in I_i$ such that
the restrictions of  $f_{\mu_i}$ and $f_{\mu_{i+1}}$ agree on
$I_i - \{a\}$.
\end{itemize}
\end{theorem}

\begin{proof} There are only finitely many points $\mu_1$ of Atiyah-Bott type
for $c$ with $\mu\ge \mu_1\ge \mu'$. Thus, it suffices to show that if
$f_\mu>f_{\mu'}$, and if
there is no point $\mu''$ of Atiyah-Bott type for $c$ with $f_\mu>f_{\mu''}
>f_{\mu'}$,
then a finite sequence of the moves listed can be applied to
$(\mu,I)$,
each one not increasing $\mu$, to
produce $(\mu',I')$. First suppose that $I'\not\subseteq I$. Then
choose  $b_0\in
I' -  I$ and define a function $f_0$ which is harmonic on
$\Delta\spcheck -  (I\cup\{b_0\})$ with $f_0(a)=f_\mu(a)$ for all $a\in I$
and $f_0(b_0)=f_{\mu'}(b_0)$. By Lemma~\ref{exten} there is a unique such
function.
 Clearly, by Proposition~\ref{abpoint}, since $\mu$ and $\mu'$ are of
Atiyah-Bott type for $c$,  the resulting function $f_0$ is determined by a
point
$\nu_0$ in $\frak h_\Ar$ such that $(\nu_0,I\cup\{b_0\})$
of Atiyah-Bott type for $c$.

Since $f_0(b_0)\le f_\mu(b_0)$, it also follows from Lemma~\ref{exten}  that
$f_\mu\ge f_0$. Since $\mu$ is harmonic at $b_0$, it follows easily
from the definitions, Lemma~\ref{inequal}, and the
fact that
$f_\mu\ge f_0$ that
$f_0$ is subharmonic at $b_0$. Suppose that  $f_0\not\ge f_{\mu'}$.
Then there exists a $b\in\Delta\spcheck$ such that
$f_0(b)< f_{\mu'}(b)$. Since
$f_0$ and
$f_{\mu'}$ are both harmonic except at $I\cup I'$, by Lemma~\ref{exten} there
exists a
$b_1\in I\cup I'$ with $f_0(b_1)< f_{\mu'}(b_1)$. Since $f_{\mu'}| I\cup
\{b_0\} \le f_0| I\cup\{b_0\}$, it follows that $b_1\in I'-(I\cup \{b_0\})$.
Now perform the same construction with $b_1$ replacing $b_0$, producing a
function
$f_1\le f_\mu$ which is harmonic except at
$I\cup \{b_1\}$ and subharmonic at $b_1$.
Let $\nu_1\in \frak h$ be the point of Atiyah-Bott type for $c$ for
which $f_{\nu_1}=f_1$.
 Since $f_1| I\cup \{b_1\}
> f_0 | I\cup\{b_1\}$ and since $f_1$ is subharmonic except at $I$ and
harmonic except at $I\cup\{b_1\}$, it follows from Lemma~\ref{inequal} that
$f_\mu \geq f_1> f_0$. Continuing in this way, we can eventually choose
$b\in I' -  I$ and a corresponding point $\nu\in \frak h_\Ar$ with
$(\nu,I\cup \{b\})$
of Atiyah-Bott type for $c$ such that the corresponding function $f_\nu$
has the property that
 $f_\mu\ge f_\nu\ge f_{\mu'}$.
Clearly with this choice, $(\nu, I\cup\{b\})$ is obtained from
$(\mu,I)$ by a move of Type (1) and  $f_\mu\ge f_\nu\ge f_{\mu'}$.
It follows from our hypothesis that either $\nu=\mu$ or
$\nu=\mu'$.
Thus,
either this move replaces $\mu$ by $\mu'$ or it leaves $\mu$ unchanged
and decreases the cardinality of $I'-I$.
We can continue applying moves of Type (1)
in this manner until $I'\subseteq I$.

Now suppose that $I'\subseteq I$. If for some $b\in I'$ we have
$f_\mu(b)>f_{\mu'}(b)$, define the function $f$ by setting
$f(a)=f_\mu(a)$ for all $a\in I - \{b\}$, $f(b)=f_{\mu'}(b)$, and requiring
that
$f$ be harmonic elsewhere. Again it is clear by Proposition~\ref{abpoint}
that $f=f_\nu$ for a point
$\nu\in \frak h_\Ar$ of Atiyah-Bott type for $c$. Since $f_\mu$ and $f_\nu$
are harmonic except at $I$, by Lemma~\ref{exten}
$f_\mu>f_\nu$. Since $I'\subseteq I$,  both $f_\nu$ and $f_{\mu'}$ are
harmonic except  at
$I$. Since $f_\mu\ge f_{\mu'}$, it follows that $f_\nu|I\ge f_{\mu'}|I$.
Applying Lemma~\ref{exten} once again, we see that $f_\nu\ge f_{\mu'}$. We can
then obtain
$(\nu, I)$ from $(\mu, I)$ by by a move of Type (3). Since $f_\mu>
f_{\nu}\ge f_{\mu'}$, it follows from our hypothesis that $\nu=\mu'$.
Thus, $(\mu',I)$ is obtained from $(\mu,I)$ by a single move of Type
(3). Thus, this allows us to arrange that $I'\subseteq I$ and that
$f_\mu|I'= f_{\mu'}|I'$.
Then, we can
obtain $\mu'$ from $\mu$ by a single move of Type (2).
\end{proof}

\subsection{Bundle deformations}

In this subsection, we shall concentrate on the case $g(C) =1$.
Here is the result which allows us to cover the  moves of
Theorem~\ref{decrease} by bundle moves.

\begin{theorem}\label{cover} Suppose that $g(C) =1$.
Fix $c\in \pi_1(G)$. Suppose
$(\mu,I)$ and $(\mu',I')$ are
of Atiyah-Bott type for some $c$, $f_\mu\ge f_{\mu'}$, and that
$(\mu',I')$ is  obtained
from $(\mu,I)$ by one of the three  moves described in
Theorem~\ref{decrease}.
Let $\Xi$ be a holomorphic  $P^I$-bundle over $C$ such that
$\eta=\Xi/U^I$ is a semistable  $L^I$-bundle whose Atiyah-Bott point is
$\mu$. Then there is an arbitrarily  small deformation of
$\Xi\times_{P^I}G$ to a $G$-bundle $\Xi'$ which has a reduction $\Xi_{P^{I'}}$
over
 $P^{I'}$ such that $\eta'=\Xi_{P^{I'}}/U^{I'}$ is a semistable
$L'$-bundle with Atiyah-Bott point $\mu'$.
\end{theorem}

\begin{proof}
We shall consider each type of move separately.
We begin with two general lemmas.

\begin{lemma}\label{paracover} Let $H_1$ and
$H_2$ be connected linear algebraic groups over $\Cee$. Suppose that $p\colon
H_1\to H_2$ is a surjective homomorphism. Let $\Xi_1$ be a holomorphic
$H_1$-bundle over a curve $C$ of arbitrary genus and let
$\Xi_2$ be the holomorphic $H_2$ bundle $\Xi_1\times_{H_1}H_2$. Then any
small deformation of $\Xi_2$ can be covered by a small deformation of
$\Xi_1$.
\end{lemma}

\begin{proof} The tangent space to the deformations of $\Xi_1$ is given by
$H^1(C;\ad_{H_1}\Xi_1)$. The tangent space to  the space of deformations
of $\Xi_2$ is given by $H^1(C;\ad_{H_2}\Xi_2)$. The natural
map of vector bundles $\ad_{H_1}\Xi \to \ad_{H_2}\Xi$ is
surjective, and hence the map $H^1(C;\ad_{H_1}\Xi)\to H^1(C;\ad_{H_2}\Xi_2)$
is surjective.  From this the result follows.
\end{proof}

The following lemma will be used for holomorphic bundles, but of course a
similar result holds in the $C^\infty$ category.

\begin{lemma}\label{fibprod} Suppose that $G_1, G_2, H$ are complex Lie groups
and that $\phi_i\colon G_i\to H$ are homomorphisms with $\phi_1$ surjective.
Let
$\zeta$ be a holomorphic principal $G_1$-bundle over a complex manifold $X$,
and suppose that $\zeta\times_{G_1}H$ is isomorphic to a bundle of the form
$\zeta'\times_{G_2}H$, where $\zeta'$ is a holomorphic $G_2$-bundle. Then
there exists a holomorphic $G_1\times_HG_2$-bundle $\zeta''$, such that
$\zeta$ is isomorphic to $\zeta''\times _{(G_1\times_HG_2)}G_1$.
\end{lemma}
\begin{proof} Choose an open cover $\{\mathcal{U}_i\}$ of $X$ and
trivializations of
$\zeta$ and
$\zeta'$  over each $\mathcal{U}_i$. Suppose that $g_{ij}$ are the transition
functions for $\zeta$ and $g_{ij}'$ are those for $\zeta'$. There exists a
$0$-cochain $\{h_i\}$ with values in $H$ such that $h_i\phi_1(g_{ij})h_j^{-1}
= \phi_2(g'_{ij})$. After shrinking the open cover, we can assume since
$\phi_1$ is surjective that the functions $h_i$ lift to functions $\widetilde
h_i\in G_1$ with $\phi_1(\widetilde h_i) = h_i$. Using the $0$-cochain
$\{\widetilde h_i\}$ to modify $\{g_{ij}\}$, we may then assume that the
$g_{ij}$ satisfy:
$\phi_1(g_{ij}) = \phi_2(g'_{ij})$. Thus, $(g_{ij}, g_{ij}') \in
G_1\times_HG_2$. Clearly, $\{(g_{ij}, g_{ij}')\}$ is a $1$-cocycle, and the
bundle $\zeta''$ which it defines satisfies: $\zeta''\times
_{(G_1\times_HG_2)}G_1 \cong \zeta$.
\end{proof}

Returning to the proof of Theorem~\ref{cover}, we begin with moves of
Type (2), the simplest to describe.
This case relies on the following elementary lemma which is  implicit in
\cite{Ra} and \cite{AtBo}:

\begin{lemma}\label{deftoss}
Let $L$ be a reductive group.
Every holomorphic $L$-bundle over $C$ has an arbitrarily
small deformation  to a semi\-stable  $L$-bundle.
\end{lemma}

\begin{proof} Fix a $C^\infty$ $L$-bundle $\eta$ over $C$. The
  $(0,1)$-connections on
$\eta$ which determine a semistable  holomorphic $L$ bundle structure are a
non-empty open subset of the space of all $(0,1)$-connections on $\eta$ and
hence form a dense open subset.
\end{proof}

For a move of Type (2)  $I'\subseteq I$  and $f_{\mu'}|I'=f_{\mu}|I'$.
Thus, $L^I\subseteq L^{I'}$ and $P^I\subseteq P^{I'}$. By
Lemma~\ref{deftoss}  the
$L^{I'}$-bundle $\eta\times_{L^I}L^{I'}$ has an arbitrarily small
deformation to a semi\-stable $L^{I'}$-bundle $\eta'$. The Atiyah-Bott
point of $\eta'$ is $\mu'$. Now apply Lemma~\ref{paracover} to the surjection
$P^{I'} \to L^{I'}$ and the bundle
$\Xi\times_{P^I}P^{I'}$ to produce an arbitrarily small deformation of the
$P^{I'}$-bundle
$\Xi\times_{P^I}P^{I'}$ to a
$P^{I'}$-bundle $\Xi'$ with $\Xi'/U^{I'}$ isomorphic to $\eta'$. Viewing this
deformation of $P^{I'}$-bundles as giving a deformation of $G$-bundles
exhibits the required deformation for a move of Type (2).

Now we turn to a move of Type (1). This time $I'=I\cup\{a\}$,
so that $L^{I'}\subseteq L^{I}$ and $P^{I'}\subseteq P^{I}$. Let $P\subseteq
L^{I}$ be the maximal parabolic $P^{I'}\cap L^I$. Its Levi factor is
$L^{I'}$. Denote by $U$ its unipotent radical.
First let us consider the case where $f_{\mu'}$ is strictly
subharmonic at $a$. This means that if $\eta'$ is a semi\-stable
$L^{I'}$-bundle with Atiyah-Bott point $\mu'$, then the
Harder-Narasimhan parabolic for the $L^I$-bundle
$\eta'\times_{L^{I'}}L^I$ is $P_-$, the opposite parabolic to $P$.
The relevant lemma for this case is the following.

\begin{lemma}\label{subharmdef} Suppose that $g(C) =1$. Let $M$ be a reductive
group and let
  $\Upsilon$ be
a semistable  $M$-bundle over   $C$.
Let $Q\subseteq M$ be a parabolic
subgroup with Levi
factor $M_1$ and unipotent radical $V$. Let $Q_-$ be the opposite parabolic
in $M$, and let $V_-$ be its unipotent radical. Suppose that $\tau$ is a
semistable $M_1$-bundle over  $C$
such that the Harder-Narasimhan parabolic of $\tau\times_{M_1}M$ is $Q_-$, and
such that $\tau\times _{M_1}M$ and $\Upsilon$ are $C^\infty$
isomorphic. Then there is
an arbitrarily small deformation of
$\Upsilon$ to a bundle of the form
$\Upsilon_Q\times_QM$ where $\Upsilon_Q$ is a holomorphic $Q$-bundle
such that $\Upsilon_Q/V$ is semi\-stable and is
$C^\infty$-isomorphic to  $\tau$.
\end{lemma}

\begin{proof} Let $\frak v$, $\frak q$, $\frak m$, $\frak v_-$ be the Lie
algebras of
$V$,  $Q$, $M$, and
  $V_-$ respectively. The direct sum decomposition $\frak m =
\frak q\oplus
\frak v_-$ is preserved by the action of $M_1$. Thus
$$\ad_M (\tau\times _{M_1}M) \cong \ad_Q(\tau\times _{M_1}Q) \oplus \frak
v_-(\tau).$$
Since $Q_-$ is the Harder-Narasimhan parabolic for $\tau\times_{M_1}M$, the
bundle $\frak v_-(\tau)$ is a direct sum of semistable bundles of positive
degrees. Since $g(C)=1$, it follows  from
stability and Serre duality that
$H^1(C; \frak v_-(\tau)) = 0$. Thus the natural map
$$H^1(C; \ad_Q(\tau\times _{M_1}Q)) \to H^1(C; \ad_M (\tau\times _{M_1}M))$$
is surjective. Since $C$ is a curve, all deformations are unobstructed, and
the map from the deformation space of the $Q$-bundle
$\tau\times _{M_1}Q$ to that of the $M$-bundle
$\tau\times _{M_1}M$ is a submersion. Thus every arbitrarily small deformation
of the $M$-bundle
$\tau\times _{M_1}M$ arises from an arbitrarily small  deformation of the
$Q$-bundle
$\tau\times _{M_1}Q$. In particular, there is an arbitrarily small
deformation of the $Q$-bundle
$\tau\times _{M_1}Q$ whose associated $M$-bundle
is semistable.

The set of all semistable $L$-bundles $C^\infty$ isomorphic to $\tau$ may be
parametrized by an irreducible scheme, in the sense that there exists
an irreducible scheme $S_0$ and an $L$-bundle over $S_0\times C$ whose
restriction to every slice $\{s\}\times C$ is semistable, and moreover
such that every semistable $L$-bundle $C^\infty$ isomorphic to $\tau$
arises in this way. By the results of
\cite[Appendix]{FMII}, since by Riemann-Roch $\dim H^1(C; \frak v(\tau))$ is
independent of $\tau$, there is an irreducible scheme
$S_1$, fibered in affine spaces over
$S_0$, which parametrizes all $Q$-bundle deformations of $\tau
\times_{M_1}Q$. By the above, there is a nonempty open subset $S_1^{ss}$ of
$S_1$ and a dominant morphism from $S_1^{ss}$ to the moduli space of all
semistable $M$-bundles of the same topological type as $\Upsilon$, and this
moduli space is irreducible.  From this, the result follows.
\end{proof}

Applying this lemma with $M=L^I$, $Q=P$, $\Upsilon=\eta$, we see that
there is an arbitrarily small deformation of the $L^I$-bundle $\eta$
to a bundle of the form $\eta_P\times_PL^I$ where $\eta_P$ is a
$P$-bundle and $\eta_P/U$ is a semi\-stable $L^{I'}$-bundle with Atiyah-Bott
point
$\mu'$.

Using  Lemma~\ref{paracover} for the surjection $P^I \to L^I$, we produce an
arbitrarily small deformation of $\Xi$ to a bundle $\Xi'$ whose reduction
modulo $U^I$ is isomorphic to $\eta_P\times_PL^I$. Since $P^{I'}\subseteq
P^I$ is the preimage of $P$ under the natural projection $P^I\to L^I$, and
hence $P^{I'} = P^I\times _{L^I}P$, it follows from Lemma~\ref{fibprod} that
$\Xi'$ can be written as
$\Xi_{P^{I'}}\times_{P^{I'}}P^I$ for some
$P^{I'}$-bundle $\Xi_{P^{I'}}$ whose reduction modulo $U^{I'}$ is
$\eta_P$. This completes the discussion of the move of Type (1) in the
case when $f_{\mu'}$ is strictly subharmonic at $a$.

Now assume that $f_{\mu'}$ is harmonic at $a$. This means that
$\mu=\mu'$.

\begin{lemma}\label{harmdef} Suppose that $g(C) =1$. Let $M$ be a reductive
group and let
  $\Upsilon$ be
a semistable  $M$-bundle over   $C$. Let
$Q\subseteq M$ be a parabolic subgroup
with Levi factor $M_1$ and unipotent radical $V$.
 Then there is
an arbitrarily small deformation of
$\Upsilon$ to a bundle of the form
$\Upsilon'\times_{M_1}M$ where $\Upsilon'$ is a semi\-stable
$M_1$-bundle whose  Atiyah-Bott point is equal to that of $\Upsilon$
under the inclusion of the center of $M$ into the center of $M_1$.
\end{lemma}

\begin{proof}
Let $\frak v$ be the Lie algebra of $V$ and $\frak v_-$ the Lie
algebra of the opposite unipotent radical.
Let $Z(M_1)$ denote the identity component of the center of $M_1$.
The $M_1$-module $\frak v$ decomposes as a direct sum
$\bigoplus_\chi{\frak v}_\chi$ where $\chi\in\widehat{Z(M_1)}$ are
  characters vanishing on the identity component of the center of
  $M$. No $\chi$ is trivial, since $\frak v$ is a direct sum of root spaces
$\frak m^\alpha$ corresponding to roots which are not trivial on $Z(M_1)$.
Since
$Z(M_1)$ is the center of
$M_1$, there is an action of the group of holomorphic $Z(M_1)$-bundles on the
space of holomorphic
$M_1$-bundles. We denote this action by
$(\lambda,\Gamma)\mapsto \Gamma\otimes \lambda$.
If $c_1(\lambda)=0$ and if $\Gamma$ is semi\-stable, then
$\Gamma\otimes\lambda$
is a semi\-stable holomorphic $M_1$-bundle  which is
$C^\infty$-isomorphic to $\Gamma$.
Fix $\Gamma$ a semi\-stable $M_1$-bundle whose Atiyah-Bott point is
that of $\Upsilon$ under the inclusion $Z(M)\subseteq Z(M_1)$.
For any $\chi\in \widehat{Z(M_1)}$ trivial on $Z(M)$,
consider the vector bundle ${\frak v}_\chi(\Gamma)$. This is a semi\-stable
vector bundle of degree zero over the genus one curve $C$.
Of course, ${\frak v}_\chi(\Gamma\otimes\lambda)={\frak
  v}_\chi(\Gamma)\otimes \chi(\lambda)$. For generic choices of
$\lambda$ with $c_1(\lambda)=0$, the cohomology group $H^1(C;{\frak
  v}_\chi(\Gamma\otimes \lambda))$ vanishes for all characters $\chi$, and
hence for such generic $\lambda$ we have $H^1(C;{\frak
  v}(\Gamma\otimes\lambda))=0$.
By duality, $H^1(C;{\frak v}_-(\Gamma\otimes\lambda))=0$ for generic
such $\lambda$. For such $\lambda$,
$$H^1(C; \ad_{M_1}(\Gamma\otimes\lambda)) = H^1(C; \ad
_M((\Gamma\otimes\lambda)\times _{M_1}M)).$$ Consider the set of all pairs
$(\Gamma,\lambda)$ for which $H^1(C;{\frak
  v}(\Gamma\otimes\lambda))=0$. As in the proof of the
previous lemma, there is an irreducible  scheme $S$ parametrizing (possibly
many-to-one) the isomorphism classes of such pairs and a dominant map from
$S$ to the moduli space of semi\-stable $M$-bundles $C^\infty$-isomorphic to
$\Upsilon$. Hence, the generic such $M$-bundle can be written as
$\Upsilon'\times_{M_1}M$ with $\Upsilon'$ a semi\-stable $M_1$-bundle
with the same Atiyah-Bott point as $\Upsilon$.
\end{proof}

Applying this lemma with $M=L^I$, $Q=P$, $\Upsilon=\eta$, we see that
there is an arbitrarily small deformation of the $L^I$-bundle $\eta$
to a bundle of the form $\eta_P\times_PL^I$ where $\eta_P/U$ is a
semi\-stable $L^{I'}$-bundle with Atiyah-Bott point $\mu$. Applying
Lemma~\ref{paracover} and  Lemma~\ref{fibprod}  as before completes the
discussion in this case.

Notice that we did not use the hypothesis that there were no
Atiyah-Bott points of type $c$ strictly between $\mu$ and $\mu'$.
We include this hypothesis in order to handle the case of genus
greater than one below.

Now we consider a move of Type (3).
The relevant result is the following.

\begin{theorem}\label{reduceby1} Let $M$ be a reductive group and let $\alpha$
be a simple root of
$M$. Let
$Q=Q^\alpha$ be the corresponding maximal parabolic, $V$ its unipotent
radical and $M_1$ its Levi factor.
  Suppose that $\eta$ is  a semistable $M_1$-bundle with Atiyah-Bott point
$\mu$ with $f_\mu(\alpha\spcheck)=q$. Then there is an arbitrarily small
$M$-deformation of
$\eta\times_{M_1}M$ to a $M$-bundle of the form $\eta_Q\times_QM$ where
$\eta_Q$ is a $Q$-bundle with
$\eta_Q/V$ a semistable  $M_1$-bundle whose Atiyah-Bott point $\mu'$ satisfies
$f_{\mu'}(\alpha\spcheck) =q-1$.
\end{theorem}

The proof of this theorem involves some new ideas and we postpone it
to the next section.
We show how this result implies Theorem~\ref{cover} for moves of Type (3).
Set $M=L^{I-\{a\}}$ and $P=P^{I-\{a\}}$ and let $U$ be the unipotent radical
of $P$. We let $Q= P^I\cap M$ and $M_1 = L^I$. Since $f_\mu(a)-f_{\mu'}(a)$ is
a positive integer, applying Theorem~\ref{reduceby1} repeatedly produces an
arbitrarily small deformation of the $M$-bundle $\eta\times_{L^{I}}M$
to a bundle of the form $\eta_Q\times_{Q}M$ where $\eta_Q$ is a
semi\-stable $Q$-bundle such that the   Atiyah-Bott point of $\eta'=\eta_Q/V$
is $\mu'$. Applying Lemma~\ref{paracover} produces an arbitrarily small
deformation of the $P$-bundle $\Xi\times_{P^I}P$ to a $P$-bundle of the form
$\Xi_P$  where $\Xi_P/U$ is isomorphic to
$\eta_Q\times_{Q}M$. It follows as before from Lemma~\ref{fibprod} that
$\Xi_P$ reduces to a
$P^I$-bundle $\Xi_{P^I}$ with $\Xi_{P^I}/U^I$ isomorphic to $\eta'$.
This concludes the proof of Theorem~\ref{cover}.
\end{proof}

\noindent{\bf Proof of Theorem~\ref{mainthm}.}
Theorem~\ref{mainthm} is an immediate consequence of Theorem~\ref{cover}
and Theorem~\ref{decrease}. \endproof

\subsection{The case of higher genus}\label{highergenus}

Here is the (weaker)  statement  in the case of higher genus.

\begin{theorem}\label{cover>1} Suppose that $g(C) \ge 2$. Let
$(\mu,I)$ and $(\mu',I')$ be
of Atiyah-Bott type for $c$, $f_\mu\ge f_{\mu'}$, and suppose that
$(\mu',I')$ is  obtained
from $(\mu,I)$ by one of the three elementary moves described in
Theorem~\ref{decrease}.
Then there exists  a holomorphic  $P^I$-bundle $\Xi$ such that
$\eta=\Xi/U^I$ is a semistable  $L^I$-bundle whose Atiyah-Bott point is
$\mu$ and an arbitrarily  small deformation of
$\Xi\times_{P^I}G$ to a $G$-bundle $\Xi'$ which has a reduction $\Xi_{P^{I'}}$
over
 $P^{I'}$ with $\eta'=\Xi_{P^{I'}}/U^{I'}$ a semistable
$L'$-bundle with Atiyah-Bott point $\mu'$.
\end{theorem}

\begin{proof}
The proof given above for moves of Types (2)
and (3) applies equally well for curves of higher genus.
For the case of moves of Type (1) we further divided into the case
when $f_{\mu'}$ was strictly subharmonic at $a$ and the case when
$f_{\mu'}$ was harmonic at $a$. The proof in this case when $f_{\mu'}$
  is harmonic reduces to the following elementary fact.
Given two reductive groups $L'\subseteq L$ and a semi\-stable
$L'$-bundle $\eta'$ whose Atiyah-Bott point lies in the center of $L$,
the $L$-bundle $\eta'\times_{L'}L$ is semi\-stable.

The remaining case is where $f_{\mu'}$ is strictly subharmonic at $a$. Let
$L^I = L$.  Let $\eta_0$ be a $C^\infty$ $L$-bundle with
Atiyah-Bott point $\mu$ and such that $\eta_0\times_{L}G \cong \xi_0$.
The set $\Delta-I$ is a set of simple roots for $L$ and determines a
fundamental Weyl chamber $\ov C_0(L)$ in $\frak h$ for the Weyl
group of $L$. This chamber contains $\ov C_0$ but will be strictly
larger than it if $I\not= \emptyset$.
Since $f_\mu$ and $f_{\mu'}$ are subharmonic at $a$
and harmonic
on $\Delta\spcheck-(I\cup \{a\})$, it follows that $\mu,\mu'$ lie in
$-\ov C_0(L)$. Our hypothesis is that there is no point
$\nu$ of Atiyah-Bott type for $c$ (and the group $G$)
with $f_\mu >f_\nu >f_{\mu'}$. Hence there is no point $\nu$ of Atiyah-Bott
type for $c_1(\eta_0)$ (and the group $L$) with $f_\mu>f_\nu>f_{\mu'}$.
Of course, under the Atiyah-Bott ordering for $L$ we have
$\mu'>\mu$ since $\mu$ indexes the stratum of semi\-stable $L$-bundles.
According to Lemma~\ref{compare}, applied to the group
$L$, this implies that there is no point of Atiyah-Bott type for
$c_1(\eta_0)$ with
$\mu'>\nu>\mu$ in the Atiyah-Bott ordering.

Let $P\subseteq L$ be the maximal parabolic subgroup whose Levi factor
is $L'$, and let $U\subseteq P$ be its unipotent radical and $\frak u$
its Lie algebra. We denote by $P_-$ the opposite parabolic, by $U_-$
its unipotent radical and by ${\frak u}_-$ the  Lie algebra of this
opposite unipotent radical.
Begin with a semi\-stable $L'$-bundle $\eta'$ with Atiyah-Bott point
$\mu'$ and such that $\eta\times_{L'}L$ is $C^\infty$ isomorphic to $\eta_0$.
The Harder-Narasimhan parabolic for $\eta'$ is $P_-$. Thus the tangent
space to the stratum containing $\eta'$ is $H^1(C; \ad _{P_-}(\eta'))$ and
the normal space is $H^1(C;{\frak u}(\eta'))$.  The bundle
${\frak u}(\eta')$ is a direct sum of semi\-stable vector bundles of
negative degrees, and so $H^0(C; {\frak
  u}(\eta')) =0$. Thus, by Riemann-Roch $H^1(C;{\frak
  u}(\eta'))\not= 0$. This means that there is an arbitrarily small
$P$-deformation $\Xi_P$
of $\eta'$ such that $\Xi_P\times_ PL$  is not contained in the
stratum containing $\eta'\times_{L'}L$.
Let $\Xi=\Xi_P\times_PL$.
According to the theorem of Atiyah-Bott, this means that the
Atiyah-Bott point $\mu'>\mu(\Xi)$  in the
Atiyah-Bott ordering.
Of course, since $\mu$ is the Atiyah-Bott point of semi\-stable
$L$-bundles we have $\mu(\Xi)\ge \mu$. It now follows from the
discussion in the previous paragraph that $\mu(\Xi)=\mu$, i.e., that
$\Xi$ is a semi\-stable $L$-bundle of the topological type of $\eta$.

Thus, we have produced an $L$-bundle $\Xi$ which is semi\-stable and
whose Atiyah-Bott point is $\mu$ which reduces to a $P$-bundle $\Xi_P$
whose associated $L^{I'}$-bundle $\Xi_P/U$ is isomorphic to $\eta'$
and hence is semi\-stable with Atiyah-Bott point $\mu'$.
This completes the proof of the theorem in this last case.
\end{proof}

\noindent{\bf Proof of Theorem~\ref{mainthm2}.}
Theorem~\ref{mainthm2} is an immediate consequence of
Theorem~\ref{cover>1} and Theorem~\ref{decrease} as well as
Theorem~\ref{mainthm1}. \endproof

\subsection{An example}

Let us give an example to show that the stronger theorem which we
proved in the case of genus one does not hold for any curve of higher
genus. The example will be for the group $SL(3)$, but surely similar
examples can be constructed for any reductive group whose derived
subgroup is of rank at least two.
Fix a smooth curve $C$ of genus at least $2$.
 Consider a rank three vector bundle
with trivial determinant which is given as an extension
$$0\to \lambda\to V\to W\to 0,$$
where $\lambda$ is a line bundle of degree $2$ and $W$ is a {\bf
  stable} rank-two bundle of degree $-2$.
Note that such bundles exist on any curve of genus greater than one.
Viewing the Cartan subalgebra of $SL(3)$ as the subspace
of $\Cee^3$ of triples $(a,b,c)$ with $a+b+c=0$, with the  Weyl
chamber
defined by $a\ge b\ge c$,
the Atiyah-Bott point of $V$ is $\mu=(2,-1,-1)$. Consider the Atiyah-Bott
point $\mu'=(1,0,-1)$. Clearly, $\mu>\mu'$. The stratum
$\mathcal{C}_{\mu'}$ consists of all rank three  vector bundles with
trivial determinant and a three-step filtration whose successive
quotients are line bundles of degrees $1,0, -1$, respectively.
If there were an arbitrarily small deformation of $V$ to a bundle with
such a filtration, then by  upper semi-continuity and the
compactness of the space of line bundles of degree $-1$ on $C$, then
there would exist a non-zero map from $V$ to a line bundle of degree
$-1$. Since $\deg \lambda =2$,
the restriction of this map to $\lambda$ is trivial. Thus, there would
be an induced non-zero map from $W$ to a line bundle of degree $-1$,
contradicting the stability of $W$.
This shows that the bundle $V$ which is in the stratum $\mathcal{C}_\mu$
is not in the closure of the stratum $\mathcal{C}_{\mu'}$. Of course,
we can deform $V$ within its stratum to a bundle $V'$ which sits
in an exact sequence
$$0\to \lambda\to V'\to W'\to 0$$
with $W'$ properly semistable.
Such a bundle $V'$ has an arbitrarily small deformation to a bundle contained
in $\mathcal{C}_{\mu'}$, so that indeed the closure of $\mathcal{C}_{\mu'}$
meets $\mathcal{C}_\mu$.

\section{Elementary modifications}

\subsection{Statement of the main theorem}

The main goal in this section is to prove Theorem~\ref{reduceby1}. For
simplicity, we change notation in that theorem, so that $M$ becomes $G$,
$Q$ becomes $P$, and so on. In fact, we show the following:

\begin{theorem}\label{mainelmod} Let
  $P=P^\alpha$
be a standard maximal parabolic subgroup of $G$, with Levi
subgroup $L$. Let
$\eta$ be a semistable $L$-bundle and let $\mu=\mu(\eta)$ be such that
$f_\mu(\alpha\spcheck ) = q$. Let $P_-$ denote the opposite parabolic to
$P$ and let $U_-$ be the unipotent radical of $P_-$. Then, possibly after
replacing $\eta$ by an arbitrarily small deformation,  there exists a
$P_-$-bundle
$\xi$ such that $\xi/U_-=\eta$ and such that $\xi\times _{P_-}G \cong
\xi_P\times_PG$, where $\xi_P$ is a  $P$-bundle such that $\xi_P/U$ is a
semistable $L$-bundle whose Atiyah-Bott point $\mu'$
satisfies
$f_{\mu'}(\alpha\spcheck) = q-1$.
\end{theorem}

\noindent{\bf Proof of Theorem~\ref{reduceby1}.} It is a standard
argument (see for example \cite[Section 4]{FMII}) that there exists a
holomorphic $G$-bundle $\Upsilon$ over $\Cee \times C$ such that
$\Upsilon|\{0\}\times C\cong \eta\times_LG$ and such that, for all $t\neq 0$,
$\Upsilon|\{0\}\times C\cong \xi\times _{P_-}G \cong \xi_P\times_PG$. Thus the
$G$-bundle
$\xi_P\times_PG$ is an arbitrarily small deformation of $\eta\times
_LG$, as claimed.
\endproof

\subsection{Definition of elementary modifications}

Let $L$ be a reductive group and let $\eta$ be a principal $L$-bundle. Let
$p_0\in C$, and let $\mathcal{U}=\{\mathcal{U}_i: i=1, \dots, n\}$ be an open
cover of $C$ such that
$p_0\in \mathcal{U}_1$ and $p_0\notin \mathcal{U}_i$ for $i>1$. Fix a
trivialization of
$\eta$ with respect to the open cover $\mathcal{U}$ and let
$g_{ij}$, $i<j$, be the transition functions. By convention, this means
that there exist local trivializations $s_i$ of $\eta|\mathcal{U}_i$ and $s_j$
of $\eta|\mathcal{U}_j$,  and $s_i\circ
s_j^{-1}|\mathcal{U}_i\cap
\mathcal{U}_j$ is left multiplication by $g_{ij}$. If we are given local
sections $\sigma_i$ with $s_i(\sigma_i(x)h) = (x, h)$, then $\sigma _i =
g_{ij}^{-1}\sigma_j$. Choose a small disk
$\mathcal{U}_0\subseteq
\mathcal{U}_1$ containing $p_0$, and such that $\mathcal{U}_0\cap
\mathcal{U}_i =\emptyset$ for $i>1$. Now define a new open cover
$\mathcal{U}' = \{\mathcal{U}_0, \mathcal{U}_1-\{p_0\}, \mathcal{U}_2, \dots,
\mathcal{U}_n\}$. The transition functions for $\eta$ relative to the open
cover $\mathcal{U}'$ are equal to $g_{ij}$ for $i,j\geq 1$, and $g_{01} = 1$.
Note that $\mathcal{U}_0 \cap (\mathcal{U}_1-\{p_0\}) =
\mathcal{U}_0-\{p_0\}$ is a punctured disk.

\begin{defn} Suppose that $\varphi\colon \Cee^*\to L$ is a $1$-parameter
subgroup, in other words a homomorphism from $\Cee^*$ to $L$. With
$\mathcal{U}_0$ as above, fix an isomorphism of $\mathcal{U}_0$ to a
neighborhood of $0\in \Cee$ sending $p_0$ to $0$, and use this
isomorphism to view $\varphi$ as a function from $\mathcal{U}_0
-\{p_0\}$ to $L$. The coordinate $t$ on $\Cee$ thus defines a
coordinate, also denoted $t$, on $\mathcal{U}_0$. Define the {\sl elementary
modification\/}
$\eta_\varphi$ of $\eta$ at $p_0$ (with respect to $\varphi$ and the given
trivialization of
$\eta$) to be the
$L$-bundle given by the following transition functions $g_{ij}'$, $i<j$, with
respect to the open cover
$\mathcal{U}'$: For $(i,j)\neq (0,1)$, $g_{ij}' = g_{ij}$, and $g_{01}'
=\varphi$.
\end{defn}

The motivation for this definition is as follows. Two bundles which differ by
an elementary modification should be isomorphic away from $p_0$ and hence
birationally isomorphic. By Iwahori's theorem, every double coset in the
space $G(\Cee[[t]])\backslash G(\Cee((t)))/G(\Cee[[t]])$ is represented by a
1-parameter subgroup, and so it is natural to use these as the gluing maps
for the new bundle.

We shall not try to describe the way $\eta_\varphi$ changes for different
choices of the trivialization or the isomorphism from $\mathcal{U}_0$ to a
neighborhood of $0\in \Cee$. It is easy to see that different choices may
lead to non-isomorphic bundles. On the other hand, the topological type of
$\eta_\varphi$ is determined by that of $\eta$ and by $\varphi$. More
precisely, we have:

\begin{lemma}\label{chAB} With notation as above, $\mu(\eta_\varphi) =
\mu(\eta) +
\pi\varphi_*(1)$, where $\pi$ is the projection
from $\mathfrak{h}$ to $\mathfrak{z}_G$.
\end{lemma}
\begin{proof} It suffices to show that, for all characters $\chi$ of $L$, we
have $\deg \chi(\eta_\varphi) = \deg \chi(\eta) + n$, where $n$ is
the integer such that $\chi\circ
\varphi(t) = t^n$. But clearly $\chi(\eta_\varphi) =
\chi(\eta)_{\chi\circ\varphi}$, and we are reduced to the case of a line
bundle, i.e. $L=\Cee^*$ and $\varphi (t) = t^n$. In this case, it is easy to
see that, if $\lambda$ is the line bundle over $C$ corresponding to
$\eta$, then the line bundle corresponding to $\eta_\varphi$ is
$\lambda\otimes \scrO_C(np_0)$, and then the statement about degrees is
clear.
\end{proof}

In practice, the groups $L$ will arise as Levi subgroups of
$G$, and so a $1$-parameter subgroup of   $L$ is also a $1$-parameter
subgroup of $G$. Thus an elementary modification of an $L$-bundle $\eta$ is
also an elementary modification of the $G$-bundle $\eta\times_LG$. It will be
important to know when this construction does not change the topological type
of $\eta\times_LG$.

\begin{lemma}\label{sametoptype} Let $\xi$ be a holomorphic
$G$-bundle, and let
$\varphi$ be a $1$-parameter subgroup of $G$ which lifts to the universal
cover
$\widetilde G$. Fix an open cover and transition functions
for $\xi$, and define $\xi_\varphi$ as above. Then
$c_1(\xi_\varphi) = c_1(\xi)$, in other words the bundles $\xi_\varphi$ and
$\xi$ are isomorphic as $C^\infty$ bundles.
\end{lemma}
\begin{proof} We may assume that the cover $\mathcal{U}$ was chosen so that
each $g_{ij}$, $(i,j)\neq (0,1)$, lifts to $\widetilde g_{ij}\in \widetilde
G$. By hypothesis,
$\varphi$ lifts to $\widetilde \varphi\in \widetilde G$. Thus, we have lifted
the transition functions $g_{ij}'$ of $\xi_\varphi$ to a collection
$\widetilde g_{ij}'\in \widetilde G$. By definition, $c_1(\xi_\varphi)$ is
the coboundary of $\{\widetilde g_{ij}'\}$, viewed as an element of $H^2(C;
Z(G))$, and similarly for $c_1(\xi)$. But clearly, since there are no
triple intersections of the $\mathcal{U}_i$ involving $\mathcal{U}_0$,
we have $c_1(\xi_\varphi) = c_1(\xi)$.
\end{proof}

In fact, the condition on $\varphi$ in the lemma is both necessary and
sufficient, although we shall not use this fact.

Returning to the case of an $L$-bundle, we shall show that we may assume that
the bundle $\eta_\varphi$ is in general semistable.

\begin{lemma}\label{elmodss} With $\eta_\varphi$ as  above, there exists an
arbitrarily small deformation
$\eta_s$ of $\eta$ and of the trivializations so that $(\eta_s)_\varphi$ is
semistable.
\end{lemma}

\begin{proof} By Lemma~\ref{deftoss}, there is an arbitrarily small
deformation of
$(\eta_\varphi)_s$ which is semistable. If $g_{ij}'$ are the
transition functions
for $\eta_\varphi$ as defined above, this means that we can find transition
functions $g_{ij}'(s)$, depending on $s$ in a small disk about $0$ in $\Cee$,
such that $g_{ij}'(0) = g_{ij}'$, $0\leq i<j$, and such that the bundle whose
transition functions are $g_{ij}'(s)$ is semistable for all $s\neq 0$. Define
$g_{ij}(s) = g_{ij}'(s)$ for $(i,j) \neq (0,1)$, and $g_{01}(s) = g_{01}'(s)
\cdot \varphi^{-1}$. Clearly the functions $g_{ij}(s)$ satisfy the cocycle
condition. For $(i,j) \neq (0,1)$, $g_{ij}(0) = g_{ij}'(0) = g_{ij}$, and
$g_{01}(0) = g_{01}'(0)
\cdot \varphi^{-1} = g_{01}'\cdot \varphi^{-1}=1=g_{01}$. Thus the functions
$g_{ij}(s)$ define a deformation $\eta_s$ of $\eta$, possibly trivial, and
$(\eta_s)_\varphi=(\eta_\varphi)_s$ is
 semistable for all $s\neq 0$.
\end{proof}

\begin{corollary}\label{choice}  Let $G$ be a reductive group and let $P=
P^\alpha$ be a standard maximal parabolic subgroup of $G$, with Levi
subgroup $L$. Let
$\eta$ be a semistable $L$-bundle and let $\mu=\mu(\eta)$ be such that
$f_\mu(\alpha\spcheck ) = q$. Let $\varphi$ be the $1$-parameter subgroup of
$G$ such that $\varphi_*(1) = -\alpha\spcheck$. Possibly after replacing
$\eta$ by an arbitrarily small deformation, we may assume that:
\begin{enumerate}
\item[\rm (i)] $\eta_\varphi$ is semistable;
\item[\rm (ii)] $\eta_\varphi\times _LG$ is
$C^\infty$ isomorphic to
$\eta\times _LG$;
\item[\rm (iii)] The Atiyah-Bott point of $\eta_\varphi$ is the unique
$\mu'\in\frak z_L$ such that $f_{\mu'}(\alpha\spcheck) = f_\mu(\alpha\spcheck)
-1$.
\end{enumerate}
\end{corollary}
\begin{proof} That we can arrange (i) after an arbitrarily small
deformation of $\eta$ follows from Lemma~\ref{elmodss}. Since
$\varphi$ lifts to $\widetilde G$, (ii) follows from Lemma~\ref{sametoptype}.
Part (iii) follows from Lemma~\ref{chAB}.
\end{proof}

\subsection{Writing the elementary modification as an extension}

Our goal now is to prove:

\begin{proposition}\label{5.3.1}
 Let $\eta$, the open cover $\mathcal{U}'$, and the transition functions
$g_{ij}$ be as
  in the beginning of this section.
Fix a nonzero $X\in \frak
g^{-\alpha}$. Let
$\{u_{ij}\}$ be the
$1$-cochain with respect to the open cover $\mathcal{U}'$ with values in the
sheaf $U_-(\eta)$ defined as follows: $u_{ij} = 1$ for $(i,j)\neq (0,1)$, and
$u_{01} = \exp(t^{-1}X)$, where $t$ is the coordinate on $\mathcal{U}_0$
defined by the inclusion of $\mathcal{U}_0$ in $\Cee$. Then $\{u_{ij}\}$ is a
cocycle. Let $\xi$ be the
$P_-$-bundle defined by the $1$-cocycle $\{g_{ij}u_{ij}\}$.
Then the bundle $\xi\times_{P_-}G$ has a
reduction to $P$, such that the associated $L$-bundle is isomorphic to
$\eta_\varphi$.
\end{proposition}

\begin{proof}
 Since there are no nonempty triple intersections
involving $\mathcal{U}_0$ and $\mathcal{U}_1$, $\{u_{ij}\}$ is vacuously a
$1$-cocycle. To prove the rest of the result, we
shall show: for all $i<j$, there exists
$v_{ij}^+\colon
\mathcal{U}_i\cap \mathcal{U}_j \to U$ and, for all $i$, there exists
$u_i\colon \mathcal{U}_i\to G$ such that
$$u_ig_{ij}'v_{ij}^+u_j^{-1} = g_{ij}u_{ij}.$$
For this says that the $1$-cocycles $g_{ij}'v_{ij}^+$ and $g_{ij}u_{ij}$ are
cohomologous.
Clearly, the cocycle $g_{ij}'v_{ij}^+$
defines a
$P$-bundle $\xi_P$ whose associated $L$-bundle is $\eta_\varphi$ as required.

We can rewrite the above condition as
$g_{ij}u_i^{g_{ij}}v_{ij}^+u_j^{-1}= g_{ij}$
if $i<j$ and $i\neq 0$, or in other words $u_i^{g_{ij}}v_{ij}^+u_j^{-1} = 1$.
This condition is automatically satisfied if $u_i\in U$ for $i\geq 1$ by
setting
$v_{ij}^+ = (u_i^{g_{ij}})^{-1}u_j$. Thus for example we can take $u_i=1$
for $i\geq 1$ and set $v_{ij}^+ = 1$ for $1\leq i\leq j$. For
$(i,j)=(0,1)$ we get the single condition
$u_0{\varphi} v_{01}^+u_1^{-1}= \exp(t^{-1}X)$ for some $v_{01}^+\colon
\mathcal{U}_0\cap \mathcal{U}_1 \to U$. Assuming that we have chosen
$u_1=1$,   we seek a function
$u_0\colon \mathcal{U}_0\to G$  such that
$${\varphi}^{-1}u_0^{-1}\exp(t^{-1}X)= v_{01}^+\in U.$$
Clearly, it suffices to solve this equation in the $SL(2)$ which is a
subgroup of the universal cover of $G$ and whose Lie algebra is $\frak
g^{-\alpha}\oplus
\Cee\cdot
\alpha\spcheck
\oplus \frak g^\alpha$. In this copy of $SL(2)$, $\varphi =
\begin{pmatrix}t^{-1}&0\\0&t\end{pmatrix}$ and we may assume
that $\exp(t^{-1}X) =\begin{pmatrix}1&0\\t^{-1}&1\end{pmatrix}$. Taking
$u_0^{-1} =
\begin{pmatrix}0&1\\-1&t\end{pmatrix}$, we have the equation
$${\varphi}^{-1}u_0^{-1}\exp(t^{-1}X) =
\begin{pmatrix}t&0\\0&t^{-1}\end{pmatrix}\begin{pmatrix}0&1\\-1&t\end{pmatrix}
\begin{pmatrix}1&0\\t^{-1}&1\end{pmatrix} =
\begin{pmatrix}1&t\\0&1\end{pmatrix},$$
as required.
\end{proof}

\noindent{\bf Proof of Theorem~\ref{mainelmod}.}
Theorem~\ref{mainelmod} follows immediately from
Corollary~\ref{choice} and
Proposition~\ref{5.3.1}.

\begin{remark} A very similar proof shows that we can begin with a $P$-bundle
$\Xi$ such that $\Xi/U \cong \eta$ and deform the $G$-bundle $\Xi\times _PG$
to a bundle of the form $\Xi'\times _PG$, where $\Xi'/U \cong \eta_\varphi$.
\end{remark}

\section{The minimally unstable strata over maximal parabolics}

\subsection{Definition of the minimally unstable strata}

\begin{defn} Let $\mu$ be an Atiyah-Bott point of type $c$ which does not lie
in $\frak z_G$. We say that $\mu$   is {\sl
minimally unstable\/} if $\mu$ is minimal among all Atiyah-Bott points of
type $c$ which do not lie in $\frak z_G$. Equivalently, the
stratum $\mathcal{C}_\mu$ is minimally unstable if, for all $\mu'\neq \mu$ of
type
$c$, the closure of $\mathcal{C}_{\mu'}$ meets $\mathcal{C}_\mu$ if and only
if $\mu'\in \frak z_G$.
\end{defn}

Thus,  $\mathcal{C}_\mu$ is minimally unstable if and only if it consists of
unstable bundles,  and, for all $\xi$ lying in $\mathcal{C}_\mu$, every small
deformation of $\xi$ either lies in
$\mathcal{C}_\mu$ or is semistable.

We turn now to a detailed and explicit discussion of the minimally
unstable strata. Not suprisingly, as the next lemma shows, these are
always associated to maximal parabolic subgroups.

\begin{lemma}\label{minpar}
Let $\xi$ be a $C^\infty$-bundle over a curve $C$. Suppose that
$\mu\in\ov C_0$ is a minimally unstable point of Atiyah-Bott type for $\xi$.
Let $P(\mu)$ be the parabolic subgroup determined by
$\mu$, i.e.  $P(\mu)=P^{I(\mu)}$ where $I(\mu)$ consists of all
$\alpha\in \Delta$ with the property that $f_\mu$ is not harmonic at
$\alpha\spcheck$. Then $P(\mu)$ is a maximal parabolic.
Furthermore, $0<f_\mu(\alpha\spcheck)\le 1$.
\end{lemma}

\begin{proof} By definition, $\#I\geq 1$, and so it suffices to show that
$\#I=1$.
 Choose
$\alpha\in I$, and let $f\colon\Delta\spcheck\to \Ar$ be define by
$f(\alpha\spcheck)=f_{\mu}(\alpha\spcheck)$ and $f$ is harmonic  outside of
$\alpha$. By Lemma~\ref{inequal}, $f\le f_{\mu}$ and
$f(\alpha\spcheck)=f_{\mu}(\alpha\spcheck)>0$. Since $f$ is harmonic except at
$\alpha$ it follows from Proposition~\ref{superpos} that $f$ is
superharmonic. Clearly, by Lemma~\ref{abpoint} $f=f_\nu$ for some point
$\nu\in \ov C_0$ of Atiyah-Bott type for $\xi$ and $f_\nu\le f_{\mu}$, so
that   $\nu\le \mu$ in the Atiyah-Bott ordering. By minimality,
$\nu=\mu$. Thus $I(\mu) = I(\nu) =\{\alpha\}$.  This proves
the first statement in the lemma. To see the second, if
$f_\mu(\alpha\spcheck)>1$, then there is a unique function $f\colon \Delta
\spcheck \to \Ar$ which is harmonic except at $\alpha\spcheck$ and such that
$f(\alpha \spcheck = f_\mu(\alpha\spcheck) -1 >0$. But then $f$ is associated
to a point $\mu'$ of Atiyah-Bott type for $c$ such that $\mu'\in \ov C_0$,
$\mu'< \mu$, and $\mu' \notin \frak z_G$. This contradicts the choice of
$\mu$.
\end{proof}

Given $c$, there is exactly one unstable stratum of the type
considered in the last lemma for each  $\alpha\in \Delta$. We label the
Atiyah-Bott point for the stratum associated to $\alpha$ by $\mu_{c,\alpha}$.
In case $c=1$, we set $\mu_{1,\alpha}=\mu_\alpha$. Our next task is to
understand the partial ordering on the
$\{\mu_{c,\alpha}\}_{\alpha\in\Delta}$.

Here is the result that allows us to determine the Atiyah-Bott partial
order  on these points.

\begin{theorem}\label{order1} Given $\alpha, \beta \in \Delta$, let
$f_\alpha$ be a superharmonic function, harmonic except at $\alpha \spcheck$,
and similarly for $f_\beta$. Then
$f_\alpha\leq
f_\beta$ if and only if
$f_\alpha(\alpha\spcheck)  \leq f_{\beta}(\alpha\spcheck)$.
\end{theorem}
\begin{proof} Clearly, if $f_\alpha\leq
f_\beta$, then
$f_\alpha(\alpha\spcheck)  \leq f_{\beta}(\alpha\spcheck)$.

To prove the converse,  we use our
previous
results on harmonic and superharmonic functions. The Dynkin diagram of
$\Delta\spcheck -\{\alpha\spcheck\}$ is a union of $t$ connected
components $C_1,
\dots, C_t$. For $1\leq i\leq t$, let $\Delta _i\spcheck$ be the set of
the
vertices of
$C_i$ together with
$\{\alpha\spcheck\}$, and assume that $\beta\spcheck \in \Delta
_1\spcheck$. For
$i>1$, both
$f_{\alpha}$ and $f_{\beta}$ are harmonic on
$\Delta_i\spcheck-\{\alpha\spcheck\}$, and
$f_{\alpha}(\alpha\spcheck)  \leq f_{\beta}(\alpha\spcheck)$. It follows from
Lemma~\ref{exten} that $f_{\alpha}(\gamma\spcheck)\leq
f_{\beta}(\gamma\spcheck)$ for all $\gamma\spcheck \in
\Delta _i\spcheck, i\neq 1$.

Now consider the restrictions of $f_{\alpha}$ and $f_{\beta}$ to
$\Delta
_1\spcheck$. The function $f_{\alpha}$ is harmonic except at
$\{\alpha\spcheck\}$, and the function $f_{\beta}$ is superharmonic
on
$\Delta_1\spcheck$. Since $f_{\alpha}(\alpha\spcheck) \leq
f_{\beta}(\alpha\spcheck)$, Lemma~\ref{inequal} implies that
$f_{\alpha}(\gamma\spcheck)\leq f_{\beta}(\gamma\spcheck)$ for
all
$\gamma\spcheck \in
\Delta _1\spcheck$. Hence, for all
$\gamma\spcheck \in
\Delta\spcheck$,   $f_{\alpha}(\gamma\spcheck)\leq
f_{\beta}(\gamma\spcheck)$. This concludes the
proof of Theorem~\ref{order1}.
\end{proof}

\subsection{The simply connected case: A characterization of the
  partial order}

In this section, we assume that $G$ is simple and simply connected.
Then the strata of $G$-bundles whose Harder-Narasimhan parabolic is a
maximal parabolic and which are of minimal positive degree with
respect to the unique dominant character of this parabolic  are
indexed by the points $\mu_\alpha$, where $\mu_\alpha$ is the unique point
such that
$\beta(\mu_\alpha) =0$ for $\beta\neq
\alpha$ (i.e.,  $f_{\mu_\alpha}$ is harmonic except at
$\alpha\spcheck$)
and $f_{\mu_\alpha}(\alpha\spcheck)=1$. By Theorem~\ref{order1},
$\mu_\alpha\leq
\mu_\beta$ if and only if
$f_{\mu_\beta}(\alpha\spcheck) \geq 1$, in other words if and only if the
coefficient
$\varpi_\alpha(\mu_\beta)$ of $\alpha \spcheck$ in $\mu_\beta$, expressed as a
sum of the simple coroots, is at least $1$. The stratum
$\mathcal{C}_{\mu_\alpha}$ consists of bundles $\xi\times _{P^\alpha}G$,
where $\xi$ is a $P^\alpha$-bundle  such that  $\xi/U =\eta$ is a semistable
$L^\alpha$-bundle with  $\deg
\eta = 1$.

\begin{defn} A simple root $\alpha$ is {\sl special\/} if
\begin{enumerate}
\item[(i)] The Dynkin diagram associated to $\Delta -\{\alpha\}$ is a
union of
diagrams of type $A$;
\item[(ii)] The simple root $\alpha$ meets each component of the Dynkin
diagram
associated to $\Delta -\{\alpha\}$ at an end of the component;
\item[(iii)] The root $\alpha$ is a long root.
\end{enumerate}
If $R$ is of type $A_n$, then every simple root  is special. All other
irreducible
root systems have a unique special simple root. It corresponds to the
unique
trivalent vertex if the Dynkin diagram is of type
$D_n, n\ge 4$ or $E_n, n=6,7,8$. For  $R=C_n, n\ge 2$ or $G_2$, it is
the
long simple root.
For $R=B_n, n\ge 2$ and $F_4$ it is the  unique
long simple root which is not orthogonal to a short simple root.
\end{defn}

Let $\alpha$ be special. As in the proof of Theorem~\ref{order1},
suppose
that  the  Dynkin diagram of $\Delta\spcheck -\{\alpha\spcheck\}$ is a
union
of $t$ connected components $C_1,
\dots, C_t$, and let $\Delta _i\spcheck$ be the set of the vertices of
$C_i$ together with
$\{\alpha\spcheck\}$. Since the Dynkin diagram of $\Delta\spcheck$ has
at most one
trivalent vertex, each $\Delta _i\spcheck$ is a chain. We may uniquely
label the
coroots of
$\Delta_i\spcheck$ as $\beta _{i,1}\spcheck, \dots,
\beta_{i,{n_i}}\spcheck=\alpha\spcheck$, where
$\langle \beta _{i,k}\spcheck, \beta_{i,k+1}\spcheck\rangle \neq 0$.

\begin{lemma}\label{anotherlemma} With $\Delta _i\spcheck$ as above,
suppose that
$f\colon \Delta\spcheck \to
\Ar$ is superharmonic and is harmonic except at $\alpha\spcheck$ and at
$\beta
_{i,1}\spcheck\in
\Delta_i\spcheck$. Let
$f_j = f|\Delta _j\spcheck$. Then $f_j$ is a linear function in the
sense that
$m=f_j(\beta_{j,k+1}\spcheck) - f_j(\beta_{j,k}\spcheck)$ is constant. If
moreover
$j\neq i$, then $m= f_j(\beta_{j, 1}\spcheck) = f_j(\alpha\spcheck)/n_j$,
where $n_j =\#\Delta _j\spcheck$, and
$f_j(\beta_{j,k}\spcheck) = km$.
\end{lemma}
\begin{proof} The assumption that $\alpha$ is long implies that
$n(\beta, \alpha)
=-1$ for every simple  coroot $\beta$ which is not orthogonal to
$\alpha$.
Thus the condition that $f$ is harmonic except at $\beta _{i,1}\spcheck$
and
$\alpha\spcheck$ implies that
$$2f_j(\beta_{j,k+1}\spcheck) = f_j(\beta_{j,k}\spcheck) +
f_j(\beta_{j,k+2}\spcheck)$$ for all $k$ with $0< k < n_j-1$. Hence
$f_j(\beta_{j,k+1}\spcheck) - f_j(\beta_{j,k}\spcheck)$ is constant. If
in
addition
$j\neq i$, then
$2f_j(\beta_{j,1}\spcheck) = f_j(\beta_{j,2}\spcheck)$ and so
$f_j(\beta_{j,k+1}\spcheck) - f_j(\beta_{j,k}\spcheck)=
f_j(\beta_{j,1}\spcheck)$
for all
$k< n_j-1$. Thus
$f_j(\beta_{j,k}\spcheck) = kf_j(\beta_{j, 1}\spcheck)$. Hence, for
$j\neq i$,
the slope of $f_j$ is
$f_j(\alpha\spcheck)/n_j$, where $n_j$ is the cardinality of $\Delta
_j\spcheck$.
\end{proof}

\begin{proposition}\label{maxvalues} Suppose that $R$ is not of type
$A_n$. Let
$\alpha$ be the special root, and let $f\colon \Delta\spcheck \to \Ar$
be
superharmonic. Then
$f(\alpha\spcheck)$ is the maximum value of $f$. More precisely, $f$
increases
{\rm{(}}weakly{\rm{)}} toward
$\alpha\spcheck$, in the sense that if $\{\beta_1\spcheck, \dots, \beta
_k\spcheck=\alpha\spcheck\}\subseteq
\Delta\spcheck$ has diagram  a chain and is numbered so that
$\langle \beta_i\spcheck, \beta_{i+1}\spcheck \rangle \neq 0$, then
$f(\beta_i\spcheck)\leq f(\beta_{i+1}\spcheck)$ for $1\leq i\leq k-1$.
\end{proposition}
\begin{proof} It suffices to prove this for a superharmonic function $f$
which is
harmonic except at a single coroot $\beta\spcheck$. If $\beta\spcheck
=\alpha\spcheck$, then we have seen that $f$ is linear and increasing
toward
$\alpha\spcheck$ on each subset
$\Delta_i\spcheck$ with its natural ordering.

Now suppose that $\beta\spcheck \neq \alpha\spcheck$. Then in particular
$\beta\spcheck$ does not correspond to a trivalent vertex of the Dynkin
diagram,
so that $\Delta\spcheck -\{\beta\spcheck\}$ has at most two connected
components
$D'$ and
$D''$. Let $(\Delta')\spcheck$ be the set of all simple coroots in $D'$
together with $\beta$ and let $(\Delta'')\spcheck$  be the set of all
simple
coroots in $D''$ together with
$\beta\spcheck$. We suppose that
$\alpha\spcheck \in (\Delta'')\spcheck$. In particular, the Dynkin
diagram of
$(\Delta')\spcheck$ is a simply laced chain (possibly consisting of a
single
element).  Since $f'=f|(\Delta')\spcheck$ is harmonic except at the
endpoint
$\beta\spcheck$, it is linear and increases  toward
$\beta\spcheck$. Thus the maximum value of $f'$ is $f'(\beta\spcheck)$.

The
restriction $f''|(\Delta'')\spcheck$ is superharmonic, and is harmonic
except
$\beta\spcheck$, which corresponds to an end vertex. As before, we write
$(\Delta'')\spcheck =
\bigcup_i(\Delta_i '')\spcheck$, where each $(\Delta_i'')\spcheck$ is a
chain
containing
$\alpha\spcheck$. We may assume that $\beta\spcheck \in
(\Delta_1'')\spcheck$. By
Lemma~\ref{anotherlemma},
$f_i'' = f''|(\Delta_i'')\spcheck$ is linear, and, for $i>1$, $f_i''$
increases
up to
$\alpha\spcheck$. Moreover the slope of $f_i''$ is
$f(\alpha\spcheck)/n_i''$.

First suppose that $\alpha\spcheck$ is trivalent and meets
$\gamma_1\spcheck,
\gamma_2\spcheck,
\gamma_3\spcheck$, where $\gamma _i\spcheck\in (\Delta _i'')\spcheck$.
The
condition that
$f$ is harmonic at
$\alpha\spcheck$ says that $(f(\alpha\spcheck) - f(\gamma_1\spcheck)) +
(f(\alpha\spcheck) - f(\gamma_2\spcheck)) + (f(\alpha\spcheck) -
f(\gamma_3\spcheck)) = f(\alpha\spcheck)$.   For
$i=2,3$, $f(\alpha\spcheck) - f(\gamma_i\spcheck)$ is the slope of
$f_i''$, and
hence is at most
$f(\alpha\spcheck)/2$. It follows that the slope $f(\alpha\spcheck) -
f(\gamma_1\spcheck)$ of
$f_1''$ is nonnegative. Thus
$f_1''$ is also increasing toward $\alpha\spcheck$. Since $f'$ increases
toward
$\beta\spcheck$, which is the end vertex of $(\Delta_1'')\spcheck$, it
follows
that
$f|(\Delta')\spcheck\cup
(\Delta_1'')\spcheck$ increases toward $\alpha\spcheck$.

Next suppose that $R$ is not simply laced and that $\alpha$ is the long
simple root which meets a short root. If $\alpha$ does not correspond to
an end of the Dynkin diagram, then $(\Delta'')\spcheck =
(\Delta_1'')\spcheck\cup
(\Delta_2'')\spcheck$, with
$\beta\spcheck
\in
(\Delta_1'')\spcheck$, and $(\Delta_2'')\spcheck-\{\alpha\spcheck\} \neq
\emptyset$. Moreover
$f_2'' = f''|(\Delta_2'')\spcheck$ is linear with  slope
$f(\alpha\spcheck)/n_2''\leq f(\alpha\spcheck)/2$. If
$\gamma _i\spcheck\in (\Delta_i'')\spcheck$ are not orthogonal to
$\alpha\spcheck$, then the harmonic condition at $\alpha\spcheck$ says
that
$2f(\alpha\spcheck) = m_1f(\gamma_1\spcheck)+m_2f(\gamma_2\spcheck)$
where
$\{m_1, m_2\} =\{1,m\}$ with $m>1$. Hence
$$(m_1+m_2-2)f(\alpha\spcheck) = m_1(f(\alpha\spcheck)
-f(\gamma_1\spcheck)) +
m_2(f(\alpha\spcheck) -f(\gamma_2\spcheck)),$$
and thus
$$m_1(f(\alpha\spcheck)-f(\gamma_1\spcheck))\geq \left(\frac{m_2}{2} +
m_1
-2\right)f(\alpha\spcheck).$$ Since $m_2/2 + m_1 -2\geq 0$, $f_1''$ is
also
increasing toward $\alpha\spcheck$, and the proof concludes as in the
trivalent
case.

The remaining case is where $\alpha$ corresponds to an end vertex of the
Dynkin
diagram. In this case, if $\gamma\spcheck$ is the unique simple coroot
not
orthogonal to
$\alpha\spcheck$, the harmonic condition reads:
$$(m-2)f(\alpha\spcheck) = m(f(\alpha\spcheck) -f(\gamma\spcheck)).$$
Once again, $f_1''$ is  increasing toward $\alpha\spcheck$, and the
proof concludes as before.
\end{proof}

\begin{corollary}\label{strictmax} If $f\colon \Delta\spcheck \to \Ar$
is
superharmonic, then
$f$ attains its maximum value on the coroot dual to a special root.
Moreover, if
$\alpha$ is special, then
$f_{\varpi_\alpha\spcheck}$ has a strict maximum at $\alpha\spcheck$.
\end{corollary}

\begin{proof} In case $R$ is not of type $A_n$, the first statement
follows from
Proposition~\ref{maxvalues} and the second from
Lemma~\ref{anotherlemma}. In case
$R$ is of type $A_n$, the first statement is trivially true since every
vertex is
special, and the second again follows from Lemma~\ref{anotherlemma}.
\end{proof}

Let $\widetilde \alpha$ be the highest root of $R$ and write
\begin{align*}
\widetilde \alpha &= \sum _{\alpha \in \Delta}h_\alpha\alpha;\\
\widetilde \alpha\spcheck
&= \sum _{\alpha \in \Delta}g_\alpha\alpha\spcheck.
\end{align*}

\begin{corollary}\label{halpha} If $\alpha$ is special, then $h_\alpha =
\max\{h_\beta: \beta \in \Delta\}$.
\end{corollary}
\begin{proof}  Since $\tilde \alpha$ is the highest root, $n(\alpha,
\tilde
\alpha) = \alpha (\tilde \alpha\spcheck)\geq 0$ for all $\alpha \in
\Delta$,
and so $f_{\widetilde \alpha
\spcheck}$ is superharmonic. The result is then immediate from the
previous
corollary.
\end{proof}

It is not in general true in the non-simply laced case that, if $\alpha$
is a
special root,
then $g_\alpha =
\max\{g_\beta: \beta \in \Delta\}$. However, if equality does not hold,
then
there is a direct argument that $\max\{g_\beta: \beta \in \Delta\} =
g_\alpha
+1$.

\begin{corollary}\label{specmin}
 Suppose that $G$ is not of type $A_n$. Let $\alpha$ be
special.
Given
$\beta_1\spcheck,
\dots, \beta _k\spcheck=\alpha\spcheck\in
\Delta\spcheck$, indexed so that
$\langle \beta_i\spcheck, \beta_{i+1}\spcheck \rangle \neq 0$ for $1\leq
i\leq
k-1$, then
$\mu_{\beta_1} > \cdots > \mu_{\alpha}$.
\end{corollary}
\begin{proof} It suffices to show that $\mu_{\beta_1} > \mu_{\beta_2}$.
By
Proposition~\ref{maxvalues} applied to the superharmonic function
$f_{\mu_{\beta_1}}$, we see that $1=f_{\mu_{\beta_1}}
(\beta_1\spcheck)\leq
f_{\mu_{\beta_1}}(\beta_2\spcheck)$. Thus by Theorem~\ref{order1},
$\mu_{\beta_1}\geq \mu_{\beta_2}$.  Since $\beta_1\neq \beta_2$, in fact
$\mu_{\beta_1}> \mu_{\beta_2}$.
\end{proof}

\begin{lemma}
Suppose that $G$ is not of type $A_n$ for any $n$, and let $\alpha\in
\Delta$ be the special vertex. Then $\mu_{\alpha}$ is the unique
minimally unstable Atiyah-Bott point for the trivial $G$-bundle over $C$.
\end{lemma}

\begin{proof}
By Lemma~\ref{minpar} any point $\mu\in \ov C_0$ which  is of
Atiyah-Bott type for the trivial bundle and indexes a minimally unstable
stratum must be of the form $\mu_\alpha$ for some $\alpha\in \Delta$.
By Corollary~\ref{specmin}, only $\mu_\alpha$ for $\alpha$ special can
index a minimally unstable stratum.
\end{proof}

For root sysytems of type $A_n$, none of the $\mu_\alpha$ are comparable, and
they are all minimal elements. For root systems of type $C_n$ or $G_2$, all of
the order relations are accounted for by Corollary~\ref{specmin}. In the
remaining cases, not all of the order relations among the elements of the set
$\{\mu_{\beta}:\beta \in \Delta\}$ are accounted for by
Corollary~\ref{specmin}.  In the case of $B_n$, where the simple roots are
given by $\alpha _1= e_1-e_2, \dots, \alpha_{n-1} = e_{n-1}-e_n, \alpha_n =
e_n$, the remaining order relations are given by:  $\mu_{\alpha_n} \geq
\mu _{\alpha_k}$ if and only if  $k\geq n/2$. Likewise, in the case of
$D_n$, if the simple roots  are given by $\alpha_1 = e_1-e_2,
\dots, \alpha_{n-1} = e_{n-1}-e_n, \alpha_n = e_{n-1}+e_n$, then the
remaining order relations are given by:
$\mu_{\alpha_{n-1}} \geq \mu _{\alpha_k}$ if and only if  $k\geq n/2$,
and
similarly for
$\mu_{\alpha_n}$. Of course, it is easy to work out all such
inequalities for $E_6$, $E_7$, $E_8$, and $F_4$ as well.

The figures at the end of the paper give pictorial representations of
the relations  between the various strata corresponding to maximal
parabolics in the case of a simply connected simple group.

\begin{corollary} Let $G$ be a simply connected simple
  group.
Let $\alpha\in \Delta$ be special. Then the stratum
$\mathcal{C}_{\mu_\alpha}$ is minimally unstable. If in addition
$G$ is
not of type $A_n$ for any $n$, then $\mathcal{C}_{\mu_\alpha}$ is
absolutely
minimal in the sense that $\mathcal{C}_{\mu_\alpha}\preceq \mathcal{C}_\mu$
for all points $\mu$ of Atiyah-Bott type for the trivial bundle.
\endproof
\end{corollary}

\subsection{The non-simply connected case}

Here we assume that $G$ is simple but not necessarily simply connected.
In this case as well we wish to find all minimally unstable strata.
Rather than work out the theory on general principles,  we shall simply
specify all minimally unstable strata. In all cases, $c$ denotes a nontrivial
element of the center of the universal cover $\widetilde G$ of $G$ and we
assume as we may that $G=\widetilde G/\langle c\rangle$.

\bigskip
\noindent $\widetilde G=SL(n)$:
\medskip

In this case we identify the center of $SL(n)$ with $\Zee/n\Zee$.
Let $c$ be a central element of order $d$.
Then there are
$n/d$ minimally unstable strata. Their Atiyah-Bott points are
$\mu_{\alpha}/d$ where $\alpha$ is any vertex with
$\varpi_\alpha(c)\equiv 1/d \pmod \Zee$.

\bigskip
\noindent $G= SO(2n+1)$:
\medskip

In this case there is a unique minimally unstable stratum.
Let $\alpha_n\in\Delta$ be the unique short simple root.
The
Atiyah-Bott point $\mu$ for the minimally unstable stratum is
$\mu_{\alpha_n}/2$.

\bigskip
\noindent $\widetilde G=Sp(2n)$:
\medskip

In this case there is a unique minimally unstable stratum. Let
$\alpha_n =2e_n\in \Delta$ be the special vertex and let $\alpha_{n-1} =
e_{n-1} -e_n$ be the unique short simple root which is not orthogonal to
$\alpha_n$. Then the Atiyah-Bott point for the minimally unstable stratum is
$\mu_{\alpha_n}/2$ if $n$ is odd and $\mu_{\alpha_{n-1}}/2$ if $n$ is even.

\bigskip
\noindent $G= SO(2n)$:
\medskip

In this case there are two minimally unstable strata, interchanged by
the outer automorphism of $SO(2n)$. Let $\alpha_{n-1}$ and $\alpha_n$ be
roots corresponding to the ``ears" of the Dynkin diagram for $D_n$. The
Atiyah-Bott points for the two minimally unstable strata are
$\mu_{\alpha_{n-1}}/2$ and $\mu_{\alpha_n}/2$.

\bigskip
\noindent $\widetilde G= Spin (4n+2)$ and $c$ has  order $4$:
\medskip

In this case there is a unique minimally unstable stratum.
 Then the
  Atiyah-Bott point is $\mu_{\beta}/4$ where $\beta$ is the root
  corresponding to an ear of the Dynkin diagram   satisfying
  $\varpi_\beta(c)=1/4$. Replacing $c$ by $-c$ changes $\beta$ to the
  other ear of the diagram.

\bigskip
\noindent $\widetilde G= Spin (4n)$ and $G\neq SO(4n)$:
\medskip

In this case there is a unique minimally unstable stratum. Its
Atiyah-Bott point is $\mu_{\alpha_{n-3}}/2$ where $\alpha_{n-3}$ is
the vertex on the long arm of the Dynkin diagram next to the trivalent
vertex.

\bigskip
\noindent $G= \ad E_6$:
\medskip

In this case there is a unique minimally unstable stratum. Its
Atiyah-Bott point is $\mu_{\alpha}/3$ where $\alpha$ is the unique
root next to the trivalent vertex on one of the long arms of the
Dynkin diagram for $E_6$ with $\varpi_{\alpha}(c)=1/3$. Replacing $c$
by $-c$ replaces $\alpha$ by the corresponding root on the other long
arm of the diagram.

\bigskip
\noindent $G= \ad E_7$:
\medskip

In this case there is a unique minimally unstable stratum. Its
Atiyah-Bott point is $\mu_{\alpha}/2$ where $\alpha$ is the vertex
adjacent to the trivalent vertex on the long arm.

\begin{figure}
\centerline{\scalebox{.70}{\includegraphics{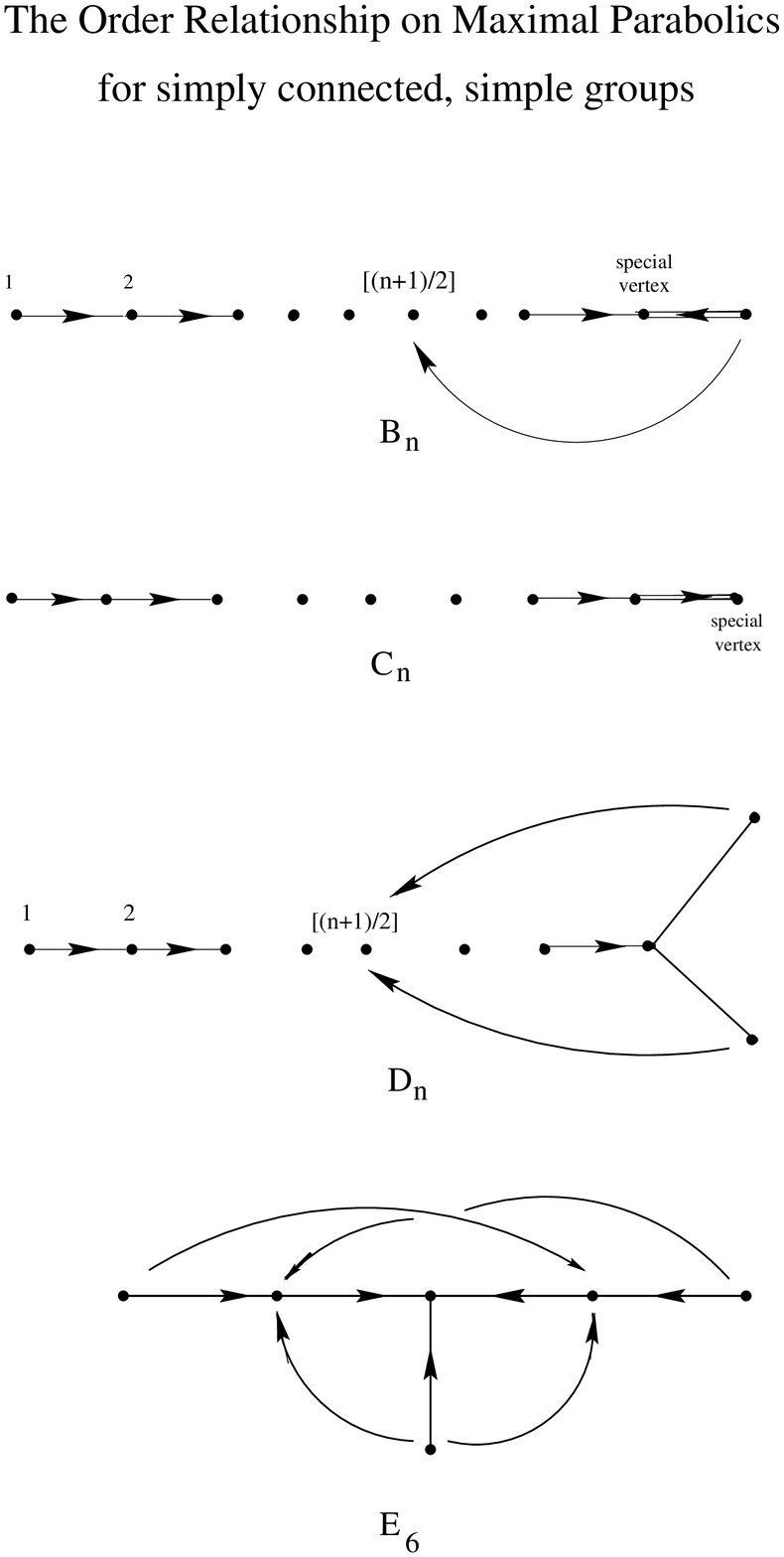}}}
\end{figure}

\begin{figure}
\centerline{\scalebox{.70}{\includegraphics{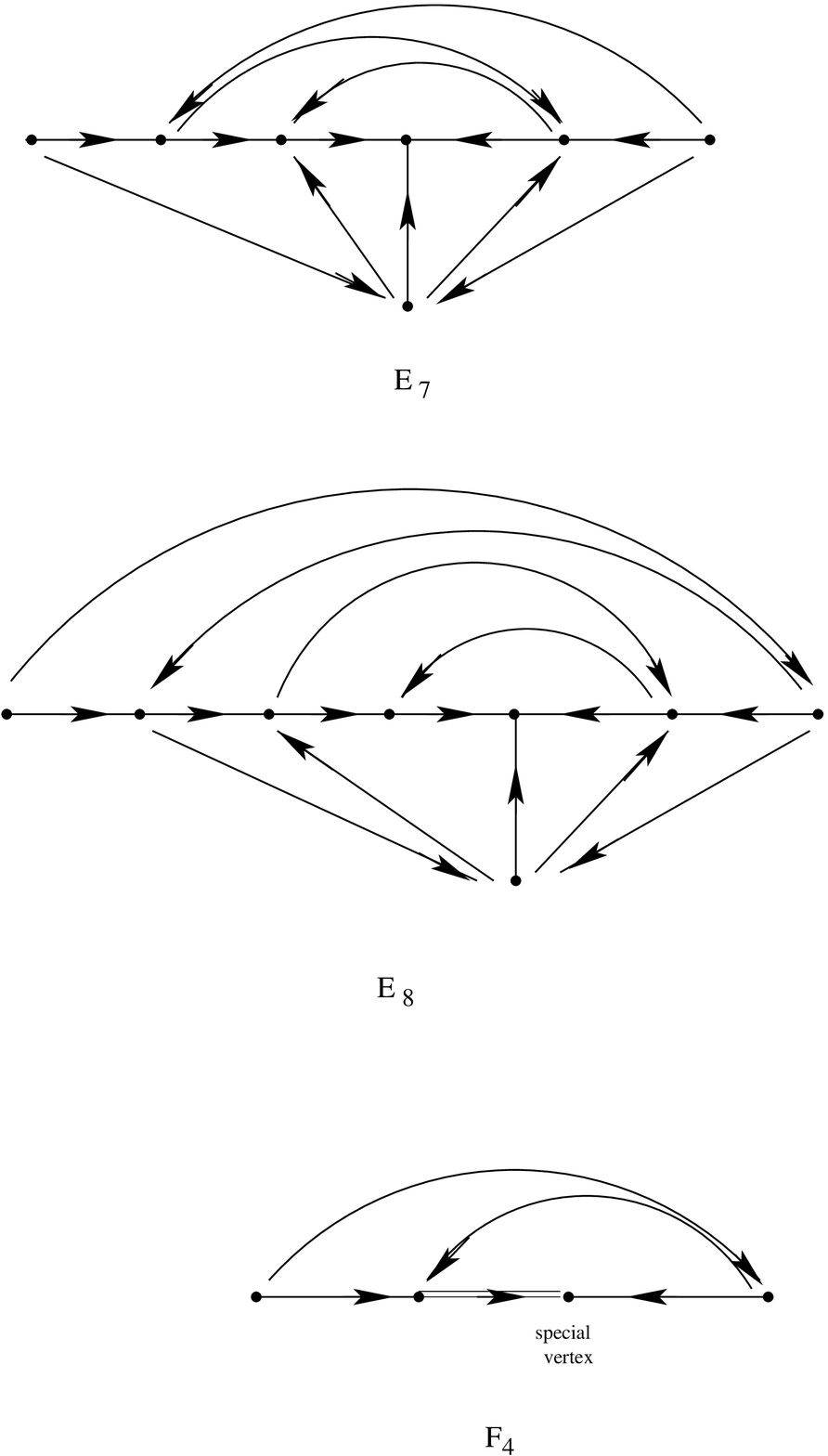}}}
\end{figure}

\bigskip
\noindent
Department of Mathematics \\
Columbia University \\
New York, NY 10027 \\
USA 

\bigskip
\noindent
{\tt rf@math.columbia.edu, jm@math.columbia.edu}


\begin{thebibliography}{10}

\bibitem{AtBo}
M. Atiyah and R. Bott,
\emph{The Yang-Mills  equations over Riemann surfaces},
Phil. Trans. Roy. Soc. London A \textbf{308}~(1982),
523--615.

\bibitem{Bour}
N. Bourbaki,
\emph{Groupes et Alg\`ebres de Lie},
Chap. 4, 5, et 6,
Masson,  Paris,
1981.

\bibitem{FM}
R. Friedman and J.W. Morgan,
\emph{Holomorphic principal bundles over elliptic curves},
AG/9811130.

\bibitem{FMII}
R. Friedman and J.W. Morgan,
\emph{Holomorphic principal bundles over elliptic curves II: The parabolic
construction}, to appear.

\bibitem{Ra}
A. Ramanathan,
\emph{Stable principal bundles on a compact Riemann surface},
 Math. Ann. \textbf{213}~(1975),
129--152.

\bibitem{Shatz}
S. Shatz,
\emph{The decomposition and specialization of algebraic families of vector
bundles},
Comp. Math. \textbf{35}~(1977),
163--187.

\end{thebibliography}
\end{document}